\documentstyle[12pt]{article}
\begin{document}

\def\tn{\otimes} \def\nc{\newcommand} \def\al{\alpha} \def\Ph{\Phi} \def\b{\beta}
\def\ph{\varphi} \def\ps{\psi} \def\inv{^{-1}} \def\ov{\over}
\def\s{\sigma} \def\iy{\infty} \def\noi{\noindent} \nc{\A}{{\rm Ai}}
\def\Dp{D^+} \def\Dm{D^-} \def\t{\tau} \def\x{\xi} \nc{\ch}{\raisebox{.3ex}{$\chi$}}
\def\pl{\partial} \def\sn{\sqrt{2n}\,} \def\D{\Delta} \def\dl{\delta} \def\r{\rho}
\renewcommand{\sp}{\vspace{1ex}} \nc{\be}{\begin{equation}} \nc{\ee}{\end{equation}}
\nc{\tl}{\tilde} \nc{\cD}{{\cal {D}}\,} \def\tr{{\rm Tr}\;}
\def\la{\lambda} \def\ph{\varphi} \def\l{\ell}
\nc{\twotwospace}[4]{\left(\begin{array}{cc}#1&#2\\&\\#3&#4\end{array}\right)}
\nc{\twotwo}[4]{\left(\begin{array}{cc}#1&#2\\#3&#4\end{array}\right)}
\nc{\twoone}[2]{\left(\begin{array}{c}#1\\\\#2\end{array}\right)}
\nc{\onetwo}[2]{\Big(#1\ \ \ #2\Big)} \def\iy{\infty} \def\cd{\cdots}
\def\dh{\widehat d} \def\hf{{1\ov 2}} \def\dg{{\rm diag}\,}
\def\e{\eta} \def\z{\zeta} \def\O{\Omega} 
\def\on{\Theta}

\hfill June 1, 2004
\begin{center}{\bf Differential Equations for Dyson Processes}\end{center}

\sp\begin{center}{{\bf Craig A.~Tracy}\\
{\it Department of Mathematics\\
University of California, Davis, CA 95616\\
email address: tracy@math.ucdavis.edu}}\end{center}
\begin{center}{{\bf Harold Widom}\\
{\it Department of Mathematics\\
University of California, Santa Cruz, CA 95064\\
email address: widom@math.ucsc.edu}}\end{center}

\begin{abstract}
We call {\it Dyson process} any process on ensembles of matrices
in which the entries undergo diffusion. We are interested in the 
distribution of the eigenvalues (or singular values) of such matrices. In the
original Dyson process it was the ensemble of $n\times n$ Hermitian matrices, and the
eigenvalues describe $n$ curves. Given sets 
$X_1,\ldots,X_m$ the probability that for each $k$ no curve passes through 
$X_k$ at time $\t_k$ is given by the Fredholm determinant of a certain matrix kernel, the 
{\it extended Hermite kernel}. For this reason we call this Dyson process the {\it Hermite process}.
Similarly, when the entries of a complex matrix 
undergo diffusion we call the evolution of its singular values
the {\it Laguerre process}, for which there is a corresponding {\it extended Laguerre kernel}.
Scaling the Hermite process at the edge leads to the {\it Airy process} (which was introduced by Pr\"ahofer 
and Spohn as the limiting stationary process for a polynuclear growth model) and 
in the bulk to the {\it sine process};
scaling the Laguerre process at the edge leads to the {\it Bessel process}.

In earlier work the authors found
a system of ordinary differential equations with independent variable $\x$
whose solution determined the probabilities
\[{\rm Pr}\,\left(A(\t_1)<\x_1+\x,\ldots,A(\t_m)<\x_m+\x\right),\]
where $\t\to A(\t)$ denotes the top curve of the Airy process.
Our first result is a generalization and strengthening of this. We assume that each $X_k$ is a 
finite union of intervals and find a system of 
partial differential equations, with the end-points of the intervals of the $X_k$ as
independent variables, whose solution determines the probability that for each $k$ no 
curve passes through $X_k$ at time $\t_k$.
Then we find the analogous systems for the Hermite process (which is more
complicated) and also for the sine process. Finally we find an analogous system of PDEs for the 
Bessel process, which is the most difficult. 
\end{abstract}

\newpage
\setcounter{equation}{0}\renewcommand{\theequation}{1.\arabic{equation}}
\begin{center}{\bf I. Introduction}\end{center}

We call {\it Dyson process\/} any process on ensembles of matrices
in which the entries  undergo diffusion. In
the  original Dyson process  \cite{D1} it was the ensemble of $n\times n$ Hermitian matrices 
$H$, where the independent coefficients of each matrix $H$   
independently execute Brownian motion
subject to a harmonic restoring force.  In one dimension this is the familiar
Ornstein-Uhlenbeck (velocity) process.  The solution to the forward (Fokker-Planck)
equation generalizes to the matrix case with the result that the probability
density of $H$
at time $\t=\t_2$  corresponding to the initial condition $H=H^\prime$ at
$\t=\t_1$ is a normalization constant depending upon $n$ and $q$ times
\[\exp\left(-{\tr(H-q H^\prime)^2\ov 1- q^2}\right), \]
where $q=e^{\t_1-\t_2}$.
As Dyson observed, the equilibrium measure as $\t_2\to\iy$ is the GUE measure of
random matrix theory.  We refer to this particular Dyson process as the
{\it Hermite process\/} for reasons that will become clear below.

With initial conditions at time $\t_1$ distributed according to the GUE
measure, the probability  that at times $\t_k$ ($k=2,\ldots, m$),\footnote
{We are assuming $\t_1<\cdots<\t_{m}$.}  $H(\t_k)$ is
in an infinitesimal neighborhood of $H_k$ 
is a normalization constant times
\be \exp\left(-\tr H_1^2\right)\,\prod_{j=2}^{m} 
\exp\left(-{\tr(H_j-q_{j-1} H_{j-1})^2\ov 1- q_{j-1}^2}\right)\, dH_1\cdots dH_{m}\,,
\label{mDistr}\ee
where $q_j=e^{\t_j-\t_{j+1}}$.  Alternatively, (\ref{mDistr}) can be interpreted
as the equilibrium measure for a chain of $m$ coupled $n\times n$ Hermitian
matrices $H_k$. 
 
In random matrix theory, and more generally Dyson processes,
one is interested in the distribution of the eigenvalues (or singular values) of $H$.  It is a classical
result of Gaudin \cite{Ga}  that the distribution functions
for the eigenvalues in GUE  are expressible in terms of the Fredholm
determinant of an integral kernel called the Hermite kernel.
In the process interpretation,  the evolution of the eigenvalues  
can be thought  of  as consisting of $n$ curves
parametrized by time. Given  $\t_1<\cd<\t_m$ and subsets $X_k$ of $\bf R$, the
quantity of interest
is the  probability that for all $k$ no curve passes through $X_k$ at time $\t_k$. 
It follows from the work of
Eynard and Mehta \cite{EM}  that this probability is also expressible
as the Fredholm determinant of
an {\it extended Hermite kernel}, an $m\times m$ matrix kernel related to the kernel associated with the 
random matrix ensemble corresponding to the equilibrium distribution.\footnote
{This was described in a lecture by Kurt Johansson \cite{jo2}, who 
recently communicated to us a sketch of his derivation \cite{jo3}. Matrix kernels,
of a different kind, also appear in \cite{AvM1}.}

Here is how it is derived. One first diagonalizes each $H_k$ and then employs the
Harish-Chandra/Itzykson-Zuber integral (see, e.g.\ \cite{ZJZ}) to integrate
out the unitary parts.  The result is that the induced measure on eigenvalues 
has a density $P(\la_{11},\ldots,\la_{1n};\ldots;\la_{m1},\ldots,\la_{mn})$ given up to
a normalization constant by
\be\prod_{k=1}^m e^{-\left({1\ov 1-q_{k-1}^2}+{q_k^2\ov 1-q_k^2}\right)
\sum_{i=1}^n \la_{k,i}^2}\, 
\prod_{k=1}^{m-1}\det\Big(e^{{2q_k\ov 1-q_k^2}\,\la_{k,i}\, \la_{k+1,j}}\Big)\,
\D(\la_1)\,\D(\la_{m}),\label{hdensity}\ee
where $q_0=q_m=0$ and $\Delta$ denotes Vandermonde determinant.\footnote
{This expression shows the connection with the theory of determinantal processes, in which
probability densities are defined by products of determinants \cite{FNH,jo1,OR}.}

In \cite{EM} it was shown that for a chain of coupled matrices with probability density
of this type the correlation functions could be expressed as block determinants whose entries
are matrix kernels evaluated at the various points, generalizing Dyson's expression for the
correlation functions for a single matrix. As with the case of random matrices, one could then
get a Fredholm determinant representation for the probability that for each $k$ 
no curve passes through $X_k$ at time $\t_k$. In the case at hand the matrix kernel 
\[ L(x,y)=(L_{ij}(x,y))_{i,j=1}^m\]
is the extended Hermite kernel and has entries
\be L_{ij}(x,y)=\left\{\begin{array}{ll}\sum\limits_{k=0}^{n-1}e^{k\,(\t_i-\t_j)}\,
\ph_k(x)\,\ph_k(y)&{\rm if}\ i\ge j,
\\&\\-\sum\limits_{k=n}^\iy e^{k\,(\t_i-\t_j)}\,\ph_k(x)\,\ph_k(y)&{\rm if}\ i<j.
\end{array}\right.\label{Lhermite}\ee
Here $\ph_k$ are the harmonic oscillator functions $e^{-x^2/2}\,p_k(x)$ where the $p_k$ are the
normalized Hermite polynomials. If
$K$ is the operator with matrix kernel $(K_{ij})$, where
\[K_{ij}(x,y)=L_{ij}(x,y)\,\ch_{X_j}(y),\]
then the probability that for each $k$ no curve passes through $X_k$ at time $\t_k$
is equal to $\det\,(I-K)$. In the special case $X_k=(\x_k,\,\iy)$ this is the probability
that the largest eigenvalue at time $\t_k$ is at most $\x_k$.

It is natural to consider also the evolution of the singular values of
complex matrices. This is the Dyson process on the space of $p\times n$ complex matrices.
(We always take $p\ge n$.) The analogue of (\ref{mDistr}) here is \cite{FNH}
\be\exp\left(-\tr A^*_1 A_1\right)\, 
\prod_{j=2}^{m} \exp\left(-{\tr\left((A_j-q_j A_{j-1})^*(A_j-q_j A_{j-1})\right)
\over 1-q_j^2}\right)\, dA_1\cdots dA_m.\label{Ldistr}\ee
After integration over the unitary parts this becomes a normalization constant times
\[\prod_{k=1}^{m}
e^{-\left({1\over 1-q_{k-1}^2}+{q_k^2\over 1-q_k^2}\right)\sum\limits_{i=1}^n\lambda_{ki}}\,
 \prod_{k=1}^{m-1}\det\left(I_\alpha\left({2 q_{k+1}\over 1-q_{k+1}^2}
\sqrt{\lambda_{k,i}\,\lambda_{k+1,j}}\right)\right)\]
\be\times\D(\lambda_{1}) \,\D(\lambda_{m}) \,
\prod_{i=1}^n \lambda_{1i}^{\alpha/2}\,
\prod_{i=1}^n \lambda_{mi}^{\alpha/2}\,\,d\lambda_{11}\cdots d\lambda_{mn},\label{Lprob}\ee
where $I_\al$ is the modified Bessel function and $\al=p-n$. (The $\la_{ki}$ are the squares
of the singular values.) This is of the same general form as for the Hermite process, and here also there is a
corresponding matrix kernel, the {\it extended Laguerre kernel}. It is given by the same
formulas (\ref{Lhermite}) as before, but now $\ph_k(x)=x^{\al/2}e^{-x/2}\,p_k(x)$
where the $p_k$ are the Laguerre polynomials $L_k^\al$, normalized.

These processes have scaling limits.
If we scale the Hermite process at the edge we obtain the 
{\it Airy process} with corresponding {\it extended Airy kernel\/}
\cite{jo1,PS} 
\be L_{ij}(x,y)=\left\{\begin{array}{ll} {\displaystyle{\int_0^\iy}} e^{-z\,(\t_i-\t_j)}\,
\A(x+z)\,\A(y+z)\,dz&{\rm if}\ i\ge j,\\&\\ {\displaystyle -\int_{-\iy}^0}e^{-z\,(\t_i-\t_j)}
\,\A(x+z)\,\A(y+z)\,dz&{\rm if}\ i<j.\end{array}\right.\label{Lairy}\ee
The Airy process consists of infinitely many curves
and as before $\det\,(I-K)$ is the probability 
that no curve passes through $X_k$ at time $\t_k$. 
In the case of greatest interest $X_k=(\x_k,\,\iy)$, and then the determinant
is equal to the probability
\be{\rm Pr}\,\left(A(\t_1)<\x_1,\ldots,A(\t_m)<\x_m\right),\label{mAiry}\ee 
where $A(\t)$ is the top curve of the Airy process. (This is what has been called 
the Airy process in the literature. It is convenient for us to use the different
terminology.)

If we scale the Hermite process in the bulk we obtain the {\it sine process} with the associated
{\it extended sine kernel}
\[ L_{ij}(x,y)=\left\{\begin{array}{ll}{\displaystyle{\int_0^1}}e^{z^2\,(\t_i-\t_j)}\,
\cos z (x-y)\,dz&{\rm if}\ i\ge j,
\\&\\-{\displaystyle{\int_1^\iy}e^{z^2\,(\t_i-\t_j)}}\,\cos z (x-y)\,dz&{\rm if}\ i<j.
\end{array}\right.\]

If we scale the Laguerre process at the bottom (the ``hard edge'') we obtain
the {\it Bessel process} and its associated {\it extended Bessel kernel}
\[ L_{ij}(x,y)=\left\{\begin{array}{ll} {\displaystyle{\int_0^1}} e^{z^2\,(\t_i-\t_j)/2}\,
\Ph_{\al}(xz)\,\Phi_{\al}(yz)\,dz&{\rm if}\ i\ge j,\\&\\ {\displaystyle -\int_1^\iy}
e^{z^2\,(\t_i-\t_j)/2}\,\Ph_{\al}(xz)\,\Ph_{\al}(yz)\,dz&{\rm if}\ i<j,
\end{array}\right.\]
where $\Ph_\al(z)=\sqrt z\,J_{\al}(z)$.\footnote
{Hints of this kernel for $m=2$ appear in \cite{Ma}.}

The Airy process $A(\t)$ was introduced by Pr\"ahofer and Spohn \cite{PS} 
as the limiting 
stationary process for a polynuclear growth model. (See also \cite{jo1}.) It is conjectured 
that it is in fact the limiting process for a wide
class of random growth models. Thus it is more significant than the Hermite
process. It might be expected that likewise the sine process (possibly) and the 
Bessel process (more likely) will prove to be more significant than the unscaled processes.

For $m=1$ the extended Airy kernel reduces to the Airy kernel and it 
is known \cite{tw3} that then (\ref{mAiry}) is expressible in terms
of a solution to Painlev\'e II.  It was thus natural  for the authors of \cite{jo1,PS}
to conjecture that the $m$-dimensional distribution functions (\ref{mAiry}) are also
expressible in terms of a solution to a system of differential equations.
This conjecture was established in two different forms, by the authors in \cite{tw22} and by 
Adler and van Moerbeke for $m=2$ in \cite{AvM}.\footnote
{In \cite{AvM1} the authors had already considered the Hermite
process in the case $m=2$, in our terminology, and 
found a PDE in $\t=\t_2-\t_1$ and the end-points of $X_1$ and $X_2$
for the probability that at time $\t_i$ no curve passes through $X_i$. In \cite{AvM} 
they deduced for the Airy process by a limiting
argument a PDE in  $\x_1,\ \x_2$ and $\t=\t_1-\t_2$ when 
$X_i=(\x_i,\,\iy)$. These equations and those we find appear to be unrelated.}
Specifically, in \cite{tw22} we found 
a system of ordinary differential equations with independent variable $\x$
whose solution determined the probabilities
\[{\rm Pr}\,\left(A_{\t_1}<\x_1+\x,\ldots,A_{\t_1}<\x_m+\x\right).\]
The $\x_k$ appeared as parameters in the equations.  

Our first result is a generalization and strengthening of this. We assume that each $X_k$ is a finite
union of intervals rather than a single interval, and find a total system of 
partial differential equations, with the end-points of the intervals of the $X_k$ as
independent variables, whose solution determines $\det\,(I-K)$. (When
$X_k=(\x_k,\,\iy)$ it it easy to recover the system of ODEs found in \cite{tw22}.)

Then we find the analogous systems for the Hermite process (which is more
complicated) and also for the sine process. Finally we find a system of PDEs for the 
Bessel process, which was the most difficult. It is possible that we could find a system 
for the Laguerre process also,
but it would be even more complicated (since Laguerre:Bessel::Hermite:Airy) and probably
of less interest.

All of these equations in a sense generalize those for the Hermite, Airy, sine and Laguerre
kernels found in \cite{tw5}, which are the cases when $m=1$. Although some of the ingredients 
are the same the equations derived here when $m=1$ are not the same as those of 
\cite{tw5}. For example, the special case of the extended Airy equations for a
semi-infinite interval and $m=1$ is the Painlev\'e~II equation whereas in \cite{tw5}
one had to do a little work to get to  Painlev\'e~II from the equations. 

We begin with Section II, where we revisit the class of probabilities for which the correlation functions
were derived in \cite{EM} and give a direct derivation of the corresponding Fredholm determinant 
representations. The method has similarities to that 
of \cite{EM} (in fact we adopt much of their notation) and the results are 
equivalent. But we avoid some awkward combinatorics. Our derivation
is analogous to that of \cite{tw14} for random matrix ensembles whereas the method of 
\cite{EM} is more like that in \cite{M}. 

In Section III we use the previous result to derive the 
extended Hermite kernel. This is of course not new. But since the derivation does not 
seem to have seen print before, this seems a reasonable place to present it. 

In the following sections we derive the systems of PDEs for the extended Airy, Hermite and
sine kernels. Presumably the other two could be obtained by scaling the equations for 
Hermite, but Airy is simpler and so we do it first. Moreover all the systems will have
the same general form, and doing Airy first will simplify the other derivations. 

In Section VII we derive the extended Laguerre kernel, and in Section VIII establish
the system of PDEs for the extended Bessel kernel.

\setcounter{equation}{0}\renewcommand{\theequation}{2.\arabic{equation}}
\begin{center}{\bf II. Extended kernels}\end{center}

For the most part we shall follow the notation in \cite{EM}. We assume the probability 
density for the eigenvalues $\la_{ki}$ ($i=1,\ldots,n,\ k=1,\ldots,m$) is given 
up to a normalization constant by
\[P(\la_{11},\ldots,\la_{1n};\ldots;\la_{m1},\ldots,\la_{mn})\]
\be=\prod_{k=1}^me^{-\sum\limits_{i=1}^nV_k(\la_{ki})}\ 
\prod_{k=1}^{m-1}\det(u_k(\la_{k,i},\,\la_{k+1,j}))\,
\D(\la_1)\,\D(\la_m),\label{density}\ee
where $V_k$ and $u_k$ are given functions satisfying some general conditions and $\D$ denotes Vandermonde
determinant. 
(Indices $i,j$ in the determinants run from 1 to $n$, and here $\la_1$ resp. $\la_m$ denotes
$\la_{1i}$ resp. $\la_{mi}$.) What we are interested in is the 
expected value of
\[\prod_{k=1}^m\prod_{i=1}^n(1+f_k(\la_{ki})),\footnote
{In our applications $f_k$ will be minus the characteristic function of $X_k$,
so the expected value will equal the probability that $\la_{ki}\not\in X_k$ for all $k$ and $i$.}
\]
so we integrate this times $P$ over all the $\la_{ki}$. 

We apply the general identity 
\[\int\cd\int\,\det(\ph_j(x_k))_{j,k=1}^n\,\cdot\,
\det(\psi_j(x_k))_{j,k=1}^n\
d\mu(x_1)\cd d\mu(x_n)\]
\[=n!\,\det\Big(\int
\ph_j(x)\,\psi_k(x)\,d\mu(x)\Big)_{j,k=1}^n\]
to the integral 
over $\la_{11},\ldots,\la_{1n}$, with the part of the
integrand containing these variables. This includes two determinants, $\D(\la_1)$
and the factor $\det\,(u_1(\la_{1i},\,\la_{2j}))$. The result is that 
this $n$-tuple
integral is replaced by the determinant
\[\det\left(\int\la_1^ie^{-V(\la_1)}\,u_1(\la_1,\,\la_{2j})\,(1+f_1(\la_1))\,d\la_1\right).\]
Then we use the same identity to rewrite the integral with respect to the $\la_{2i}$ using
this determinant and the factor $\det\,(u_2(\la_{2i},\,\la_{3j}))$. And so on. At the
end we use the determinant coming from the previous use of the identity and 
$\D(\la_m)$. The end result is that the expected value in question
is a constant times the determinant of the matrix with $i,j$ entry
\be\int\cd\int\la_1^i\la_m^j\,e^{-\sum\limits_{k=1}^mV_k(\la_k)}\,
\prod\limits_{k=1}^{m-1}u_k(\la_k,\,\la_{k+1})\prod_{k=1}^m(1+f_k(\la_k))\,d\la_1\cd
 d\la_m.\label{ev}\ee

By changing the normalization factor we may replace $\la_1^i$ by any sequence of polynomials,
which we call $P_{1i}(\la_1)$, and replace $\la_m^j$ by any sequence of polynomials, which we 
call $Q_{mj}(\la_m)$. We choose them so that after these replacements the integral with all the
$f_k$ set equal to zero equals $\dl_{ij}$. In particular the normalization constant is now equal to 1.

If we write  
\be e^{-\sum\limits_{k=1}^mV_k(\la_k)}\,\prod\limits_{k=1}^{m-1}u_k(\la_k,\,\la_{k+1})
=E_{12}(\la_1,\la_2)\,E_{23}(\la_2,\la_3)\cd E_{m-1,m}(\la_{m-1},\la_m)\label{factor}\ee
(there is some choice in  the factors on the right), we see that
the matrix in question equals the identity matrix plus the matrix with $i,j$ entry
\[\int\cd\int P_{1i}(\la_1)\prod_{k=1}^{m-1}E_{k,k+1}(\la_k,\la_{k+1})\left[\prod_{k=1}^m(1+f_k(\la_k))-1\right]
Q_{mj}(\la_m)\,d\la_1\cd d\la_m.\]

The bracketed expression may be written as a sum of products,
\[\sum_{r\ge1} \sum_{k_1<\cd<k_r}f_{k_1}(\la_1)\cd f_{k_r}(\la_{k_r}).\]
Correspondingly the integral is a sum of integrals. Consider the integral corresponding
to the above-displayed summand. For $k>j$ we define
\[E_{jk}(\la_j,\la_k)=E_{j,j+1}(\la_j,\la_{j+1})\ast\cd\ast E_{k-1,k}(\la_{k-1},\la_k),\]
where the asterisk denotes kernel composition, and set
\[P_{ki}(\la_k)=\int P_{1i}(\la_1)\,\,E_{1k}(\la_1,\la_k)\,d\la_1,\ \ \
Q_{kj}(\la_k)=\int\,E_{km}(\la_k,\la_m)\,Q_{mj}(\la_m)\,d\la_m.\]
By integrating first with respect to the $\la_k$ with $k\ne k_1,\ldots,k_r$, we see
that the corresponding integral is equal to
\[\int\cd\int f_{k_1}(\la_{k_1})\,P_{k_1,i}(\la_{k_1})\,E_{k_1,k_2}(\la_{k_1},
\la_{k_2})f_{k_2}(\la_{k_2})\cd\]
\[\cd E_{k_{r-1},k_r}(\la_{k_{r-1}},\la_{k_r})f_{k_r}(\la_{k_r})\,Q_{k_r,j}(\la_{k_r})\,d\la_{k_1}\cd d\la_{k_r}.\]

We deliberately distributed the $f$ factors as we did since if we let $A_{k,\l}$ be the operator with kernel $A_{k\l}(\la_k,\la_\l)=
E_{k\l}(\la_k,\la_\l)f(\la_\l)$ then the above may be written as the single integral
\[\int f_{k_1}(\la)\,P_{k_1,i}(\la)\,A_{k_1,k_2}\cd A_{k_{r-1},k_r}Q_{k_r,j}(\la)\,d\la.\]
(If $r=1$ we interpret the operator product to be the identity.) Replacing the index $k_1$ by 
$k$ and changing notation, we see that the sum of all of 
these equals
\[\int\sum_kf_k(\la)\,P_{k,i}(\la)\left(\sum_{r\ge0}\sum_{k_1,\ldots,k_r}
A_{k,k_1}A_{k_1,k_2}\cd A_{k_{r-1},k_r}Q_{k_r,j}(\la)\right)\,d\la,\]
where the inner sum runs over all $k_r>\cd>k_1>k$. (If $r=0$ the inner sum is interpreted to be
$Q_{k,j}(\la)$.)

We think of $f_k(\la)\,P_{k,i}(\la)$ as the $k$th entry of a row matrix and the inner sum
\[\sum_{r\ge0}\sum_{k_1,\ldots,k_r}A_{k,k_1}A_{k_1,k_2}\cd A_{k_{r-1},k_r}Q_{k_r,j}(\la)\]
as the $k$th entry of a column matrix. The integrand is the product of these matrices. If we use
the general fact that $\det\,(I+ST)=\det\,(I+TS)$ we see that the determinant of $I$ plus the matrix with the above $i,j$ entry
is equal to the determinant of $I$ plus the operator with matrix kernel having $k,\l$ entry
\[\sum_{j=0}^{n-1}\left(\sum_{r\ge0}\sum_{k_1,\ldots,k_r}A_{k,k_1}A_{k_1,k_2}\cd A_{k_{r-1},k_r}
Q_{k_r,j}(\la)\right)\,P_{\l,j}(\mu)f_\l(\mu),\]
where in the inner sum $k_r>\cd >k_1>k$.

This is the $k,\l$ entry of a certain operator matrix acting from the left on the matrix with
$k,\l$ entry
\[\sum_{j=0}^{n-1}Q_{k,j}(\la)\,P_{\l,j}(\mu)f_\l(\mu).\]
That matrix is upper-triangular, all diagonal entries are $I$, and for $k<\l$ the $k,\l$ entry
equals
\[\sum_{k<k_1<\cd <k_r<\l}A_{k,k_1}A_{k_1,k_2}\cd A_{k_r,\l}.\]
Elementary algebra shows (even for non-commuting variables $A_{k\l}$) that this is the inverse
of the upper-triangular matrix with diagonal entries $I$ and $k,\l$ entry $-A_{k\l}$ otherwise.

If we recall that $A_{k\l}(\la,\mu)=E_{k\l}(\la\,\mu)\,f_\l(\mu)$ then we see that we have shown the
following: Let $H(\la,\mu)$ be the matrix kernel given by
\[H_{k\l}(\la,\mu)=\sum_{j=0}^{n-1}Q_{k,j}(\la)\,P_{\l,j}(\mu),\]
let $E$ be the matrix kernel with $k,\l$ entry 
$E_{k\l}(\la,\mu)$ (thought of as 0 when $k\ge l$), and let $f(\mu)=\dg(f_k(\mu))$.
Then the expected value equals the determinant of
\[I+(I-Ef)^{-1}Hf=(I-Ef)^{-1}\,[I+(H-E)f].\]

The factor on the left equals $I$ plus a strictly upper-triangular matrix, so its determinant 
equals one. Therefore the expected value equals
\[\det\,[I+(H-E)f],\]
and $H-E$ is the extended kernel.
\pagebreak
\setcounter{equation}{0}\renewcommand{\theequation}{3.\arabic{equation}}
\begin{center}{\bf III. The extended Hermite kernel}\end{center}

We have times $\t_1<\cd<\t_m$ and we set 
$q_k=e^{\t_k-\t_{k+1}}$, with the conventions $\t_0=-\iy,\ \t_{m+1}=+\iy$ so that 
$q_0=q_m=0$. For the Hermite process the probability density is given by (\ref{hdensity}) so we
are in the case where
\[V_k(\la)=\left({1\ov1-q_{k-1}^2}+{q_k^2\ov 1-q_k^2}\right)\la^2,\ \ \ 
u_k(\la,\,\mu)=\exp\left\{{2q_k\ov 1-q_k^2}\la\mu\right\},\]
and we want to compute the kernel $H-E$ of Section~II.

We define the Mehler kernel
\[K(q;\,\la,\mu)=(\pi(1-q^2))^{-1/2}\,e^{-{q^2\ov 1-q^2}\la^2-{1\ov 1-q^2}\mu^2+
{2q\ov 1-q^2}\la\mu},\]
which has the representation
\be K(q;\,\la,\mu)=\sum_{i=0}^\iy q^i\,p_i(\la)\,p_i(\mu)\,e^{-\mu^2},\label{Msum}\ee
so
\be\int K(q;\,\la,\mu)\,p_i(\mu)\,d\mu=q^i\,p_i(\la).\label{Mint}\ee
Here $p_i$ are the normalized Hermite polynomials.

We can write the exponent on the left side of (\ref{factor}) as 
\[-\sum_{k=1}^{m-1}{q_k^2\ov 1-q_k^2}\la_k^2-\sum_{k=0}^{m-1}{1\ov 1-q_k^2}\la_{k+1}^2,\]
and so, aside from a normalization constant, the left side of (\ref{factor}) is equal to
\[e^{-\la_1^2}\prod_{k=1}^{m-1}K(q_k;\,\la_k,\la_{k+1}).\]
Thus we may take in (\ref{factor}) 
\[E_{12}(\la_1,\la_2)=e^{-\la_1^2}\,K(q_1;\,\la_1,\la_2),\]
\[E_{k,k+1}(\la_k,\la_{k+1})=K(q_k;\,\la_k,\la_{k+1}),\ \ \  (k>1).\]

It follows from (\ref{Msum}) that $K(q)\ast K(q')=K(qq')$
when $q,\,q'>0$ and so 
\[E_{1k}(\la,\mu)=e^{-\la^2}\,K(q_1\cd q_{k-1};\,\la,\mu).\]
In particular we deduce from (\ref{Mint}) that
\[\int\int p_i(\la)\,E_{1m}(\la,\mu)\,p_j(\mu)\,d\mu d\la=(q_1\cd q_{m-1})^j
\int p_i(\la)\,e^{-\la^2}\,p_j(\la)\,d\la=(q_1\cd q_{m-1})^j\dl_{ij}.\]
Hence we may take
\[P_{1i}=p_i,\ \ Q_{mj}=(q_1\cd q_{m-1})^{-j}\,p_j\]
as the polynomials in the previous discussion. We see that
\[P_{ki}(\mu)=\int p_i(\la)\,e^{-\la^2}\,K(q_1\cd q_{k-1};\,\la,\mu)\,d\la
=(q_1\cd q_{k-1})^i\,p_i(\mu)\,e^{-\mu^2},\ \ \ (k>1),\]
\[Q_{kj}(\la)=\int K(q_k\cd q_{m-1};\,\la,\mu)\,Q_{mj}(\mu)d\mu=
(q_1\cd q_{k-1})^{-j}\,p_j(\la),\ \ \ (k>1),\]
\[Q_{1j}(\la)=\int e^{-\la^2}K(q_1\cd q_{m-1};\,\la,\mu)\,Q_{mj}(\mu)d\mu
=e^{-\la^2}\,p_j(\la).\]
It follows that $H$ is the matrix with $k,\l$ entry
\[\sum_{j=0}^{n-1}\left({q_1\cd q_{\l-1}\ov q_1\cd q_{k-1}}\right)^j\,p_j(\la)\,p_j(\mu)\]
left-multiplied by the matrix $\dg(e^{-\la^2}\ 1\ \cd\ 1)$ and right-multiplied
by the matrix \linebreak $\dg(1\ e^{-\mu^2}\ \cd \ e^{-\mu^2})$. Similarly
$E$ is the strictly upper-triangular matrix with $k,\l$ entry 
$K(q_k\cd q_{\l-1};\,\la,\mu)$
left-multiplied by the matrix $\dg(e^{-\la^2}\ 1\ \cd\ 1)$.

Thus we have computed $H-E$. The actual extended Hermite kernel will be a modification of this.
The determinant is unchanged if we multiply $H-E$ on the left by 
$\dg(e^{\la^2/2}\ e^{-\la^2/2}\ \cd\ e^{-\la^2/2})$ and on the right by
$\dg(e^{-\mu^2/2}\ e^{\mu^2/2}\ \cd\ e^{\mu^2/2})$. Recalling that $\ph_i$ are
the harmonic oscillator functions
and recalling the definition of the $q_k$ in terms
of the $\t_k$ we see that the expected value in question is equal to 
\[\det\,[I+(\hat H-\hat E)f],\]
where 
\[\hat H_{k\l}(\la,\mu)=\sum_{j=0}^{n-1}e^{j\,(\t_k-\t_\l)}\,\ph_j(\la)\,\ph_j(\mu),\]
and $\hat E$ is the strictly upper-triangular matrix with $k,\l$ entry 
\[e^{(\mu^2-\la^2)/2}\,K(e^{\t_k-\t_\l}\,;\,\la,\mu).\]
If we observe that by (\ref{Msum})
\[e^{(\mu^2-\la^2)/2}\,K(e^{\t_k-\t_\l}\,;\,\la,\mu)=
\sum_{j=0}^\iy e^{j\,(\t_k-\t_\l)}\,\ph_j(\la)\,\ph_j(\mu)\]
when $k<\l$ we see that $\hat H-\hat E$ has $k,\l$ entry
\[(\hat H-\hat E)_{k\l}=\left\{\begin{array}{ll}
\sum\limits_{j=0}^{n-1}e^{j\,(\t_k-\t_\l)}\,\ph_j(\la)\,\ph_j(\mu)&{\rm if}\ k\ge\l,\\&\\
-\sum\limits_{j=n}^\iy e^{j\,(\t_k-\t_\l)}\,\ph_j(\la)\,\ph_j(\mu)&{\rm if}\ k<\l,\end{array}\right.\]
which is the extended Hermite kernel (\ref{Lhermite}).

\sp
\setcounter{equation}{0}\renewcommand{\theequation}{4.\arabic{equation}}

\begin{center}{\bf IV. PDEs for the extended Airy kernel}\end{center}

We consider first the case $X_k=(\x_k,\iy)$, so that
\[\det\,(I-K)={\rm Pr}\,\left(A(\t_1)<\x_1,\ldots,A(\t_m)<\x_m\right).\] 
The derivation is simplest here but
it will also give the main ideas for all the derivations.

Observe first that 
\be\pl_k\,K=-L\,\dl_k,\label{plK}\ee
where $\dl_k$ denotes multiplication by the diagonal matrix with all entries zero 
except for the $k$th, which equals $\dl(y-\x_k)$. It follows that
if we let $R=K\,(I-K)\inv$, then
\[\pl_k\log\,\det(I-K)=-{\rm Tr}\;(I-K)\inv\,\pl_kK=R_{kk}(\x_k,\,\x_k).\]
The matrix entries on the right will be among the unknowns. To explain the others,
let $A(x)$  denote the 
$m\times m$ diagonal matrix 
$\dg(\A(x))$ and $\ch(x)$ the diagonal matrix
$\dg(\ch_k(x))$, where $\ch_k=\ch_{(\x_k,\iy)}$. Then we define the matrix functions
$Q(x)$ and $\tl Q(x)$ by
\[Q=(I-K)\inv A,\ \ \ \tl Q=A\,\ch\,(I-K)\inv\]
(where for $\tl Q$ the operators act on the right).
These and $R(x,y)$ are functions of the $\x_k$ as well as $x$ and $y$. We
define the matrix functions $q,\ \tl q$ and $r$ of the $\x_j$ only by
\[q_{ij}=Q_{ij}(\x_i),\ \ \ \tl q_{ij}=\tl Q_{ij}(\x_j),\ \ \ r_{ij}=R_{ij}(\x_i,\,\x_j).
\footnote{At points of discontinuity we always take limits from the right. For example we interpret 
$R_{ij}(x,\,\x_j)$ as the limit $R_{ij}(x,\,\x_j+)$.}\]
Our unknown functions will be these and the matrix functions $q'$ and $\tl q'$ defined by
\[q'_{ij}=Q'_{ij}(\x_i),\ \ \ \tl q'_{ij}=\tl Q_{ij}'(\x_j).\]
We shall also write $r_x$ and $r_y$ for the matrices $(R_{xij}(\x_i,\x_j))$ and 
$(R_{yij}(\x_i,\x_j))$.

The $\x_k$ are the independent variables in our equations. We denote by $\x$ the
matrix $\dg(\x_k)$ and by $d\x$ the
matrix of differentials $\dg(d\x_k)$. With these notations our system of equations is
\begin{eqnarray}
dr&=&-r\,d\x\,r+d\x\,r_x+r_y\,d\x,\label{apde1}\\
dq&=&d\x\,q'-r\,d\x\,q,\label{apde2}\\
d\tl q&=&\tl q'\,d\x-\tl q\,d\x\,r,\label{apde3}\\
dq'&=&d\x\,\x\,q-(r_x\,d\x+d\x\,r_y)\,q+d\x\,r\,q',\label{apde4}\\
d\tl q'&=&\tl q\,\x\,d\x-\tl q\,(d\x\,r_y+r_x\,d\x)+\tl q'\,r\,d\x.\label{apde5}
\end{eqnarray}

One sees that the right sides involve the diagonal entries of $r_x+r_y$ and the off-diagonal entries of $r_x$ and
$r_y$. We shall show below that these are ``known'' in the sense that they are expressible 
algebraically in terms of our unknown
functions, so the above is a closed system of PDEs.

We begin by establishing the assertions about $r_x$ and $r_y$.

In the following $D=d/dx$, we set $\r=(I-K)\inv$ and 
$\dl=\sum_k\dl_k$, and $\t$ is the diagonal matrix $\dg(\t_k)$. 
We denote by $\on$ the matrix with all entries
equal to one.
For clarity we sometimes  write the kernel of an operator
in place of the operator itself.

\sp
\noi{\bf Lemma 1}. We have the commutator relation 
\be[D,\,R]=-Q(x)\,\on\,\tl Q(y)+R\,\dl\,\r+[\t,\,R].\label{DRcom}\ee
\sp
\noi{\bf Proof}. Integrating by parts in (\ref{Lairy}) gives
\[ [D,\,K]_{ij}=-\A(x)\,\A(y)\,\ch_j(y)+L_{ij}(x,\x_j)\,\dl(y-\x_j)
+(\t_i-\t_j)\,K_{ij}(x,\,y).\]
Equivalently,
\[ [D,\,K]=-A(x)\,\on\, A(y)\,\ch(y)+L\,\dl+[\t,\,K].\]
To obtain $[D,\,R]$ we replace $K$ by $K-I$ in the commutators and left- and 
right-multiply by $\r$. The result is (\ref{DRcom}).\footnote
{Because of the fact $\r\,L\,\ch=R$ and our interpretation of $R_{ij}(x,\,\x_j)$ as 
$R_{ij}(x,\,\x_j+)$
we are able to write $R\,\dl\,\r$ in place of $\r\,L\,\dl\,\r$.}
\sp

If we take the $i,\,j$ entry of both sides of (\ref{DRcom}) and set $x=\x_i,\ y=\x_j$ we obtain
\be r_x+r_y=-q\,\on\,\tl q+r^2+[\t,\,r].\label{plrsum}\ee
Thus all entries of $r_x+r_y$ are known.

For the off-diagonal entries of $r_x$ and $r_y$ we need a second commutator identity.
Here $M$ is multiplication by $x$.
\sp

\noi{\bf Lemma 2}. We have 
\[ [D^2-M,\,\r]=R\,\dl\,\r_x-R_y\,\dl\,\r,\]
where $R_y(x,y)$ is interpreted as not containing a delta-function summand.
\sp

\noi{\bf Proof}. We use the facts that $D^2-M$ commutes with $L$ and that $M$ commutes with $\ch$. 
These give
\[ [D^2-M,\,K]=[D^2-M,\,L\,\ch]=L\,[D^2-M,\,\ch]=L\,[D^2,\,\ch]=L\,(\dl\, D+D\,\dl).\]
Using the commutator identity
\[[T,\,(I-K)\inv]=(I-K)\inv\,[T,\,K]\,(I-K)\inv,\]
valid for any operators $T$ and $K$, we deduce
\[ [D^2-M,\,\r]=\r\,L\,\dl\, D\,\r+\r\,L\,D\,\dl\,\r.\]
The first term on the right equals $R\,\dl\,\r_x$. The second term equals $-R_y\,\dl\,\r$ where 
$R_y$ is interpreted as not containing the delta-function summand.
This establishes the lemma.
\sp

Lemma 1 says 
\[R_x+R_y=-Q(x)\,\on\,\tl Q(y)+R\dl\r+[\t,\,R],\]
and applying $\pl_x-\pl_y$ to both sides gives
\[R_{xx}-R_{yy}=-Q'(x)\,\on\,\tl Q(y)+Q(x)\,\on\,\tl Q'(y)+R_x\dl\r-R\dl\r_y+[\t,\,R_x-R_y].\]
Lemma 2 says
\[R_{xx}-R_{yy}-(x-y)\,R=R\dl\r_x-R_y\dl\r.\]
Equating the two expressions for $R_{xx}-R_{yy}$ gives
\be (x-y)\,R(x,y)=-Q'(x)\,\on\,\tl Q(y)+Q(x)\,\on\,\tl Q'(y)+(R_x+R_y)\dl\r-R\dl(\r_x+\r_y)
+[\t,\,R_x-R_y].\label{MRcom}\ee
Taking the $i,j$ entries and setting
$x=\x_i$, $y=\x_j$ give 
\[[\x,\,r]+r\,r_x-r_y\,r=-q'\,\on\,\tl q+q\,\on\,\tl q'+r_x\,r-r\,r_y+[\t,\,r_x-r_y],\]
or
\be[\t,\,r_x-r_y]=q'\,\on\tl q-q\,\on\,\tl q'+[r,\,r_x+r_y]+[\x,\,r].\label{rxrycom}\ee
The left side 
has $i,j$ entry $(\t_i-\t_j)\,(r_{xij}-r_{yij})$ and the right side is known.\footnote
{Here $r_{xij}$ is notational shorthand for $(r_x)_{ij}$ and $r_{yij}$ for $(r_y)_{ij}$.}
Therefore the off-diagonal entries of
$r_x-r_y$ are known, and therefore so also are
the off-diagonal entries of $r_x$ and $r_y$ individually.

To be more explicit we define matrices $U$ and $V$ by 
\[U=-q\,\on\,\tl q+r^2+[\t,\,r]\]
and 
\[V_{ij}={(q'\,\on\,\tl q-q\,\on\,\tl q')_{ij}+[r,\,-q\,\on\,\tl q+[\t,\,r]-\x]_{ij}\over \t_i-\t_j}\]
when $i\ne j$. Then (\ref{plrsum}) says
\[r_{xij}+r_{yij}=U_{ij}\]
and (\ref{rxrycom}) gives
\[r_{xij}-r_{yij}=V_{ij}\]
when $i\ne j$. It follows that for such $i,j$ we have
\[d\x_i\,r_{xij}+d\x_j\,r_{yij}={1\ov2}(d\x_i+d\x_j)U_{ij}
+{1\ov2}(d\x_i-d\x_j)V_{ij},\]
\[d\x_i\,r_{yij}+d\x_j\,r_{xij}={1\ov2}(d\x_i+d\x_j)U_{ij}
-{1\ov2}(d\x_i-d\x_j)V_{ij}.$$
The same hold when $i=j$ if we interpret the second terms to be zero then. More
succinctly,
\[d\x\,r_x+r_y\,d\x={1\ov2}\{d\x,\,U\}+{1\ov2}[d\x,\,V],\ \ \ 
d\x\,r_y+r_x\,d\x={1\ov2}\{d\x,\,U\}-{1\ov2}[d\x,\,V],\]
where the curly brackets indicate anticommutator.
These give the
explicit representations for the terms involving $r_x$ and $r_y$ in the equations.
\sp

With our assertions concerning $r_x$ and $r_y$ established we proceed to derive the equations.
It follows from the general identity
\[\pl_k\,(I-K)\inv=(I-K)\inv\,\pl_kK\,(I-K)\inv,\]
relation (\ref{plK}) and the remark in footnote 7 that
\be \pl_k\, \r=-R\,\dl_k\,\r.\label{plrho}\ee
{}From this we obtain (since $\pl_k\,R=\pl_k\,\r$) 
\[\pl_k \,r_{ij}=\pl_k\,(R_{ij}(\x_i,\x_j))=(\pl_k\,R_{ij})(\x_i,\x_j)+
R_{xij}(\x_i,\x_j)\,\dl_{ik}+R_{yij}(\x_i,\x_j)\,\dl_{jk}\]
\[=-r_{ik}\,r_{kj}+R_{xij}(\x_i,\x_j)\,\dl_{ik}+R_{yij}(\x_i,\x_j)\,\dl_{jk}.\]
Multipliying by $d\x_k$ and summing over $k$ give (\ref{apde1}). 
 
Using (\ref{plrho}) applied to $A$ we obtain
\be\pl_k \,q_{ij}=Q'_{ij}(\x_i)\,\dl_{ik}-(R\,\dl_k\,Q)_{ij}(\x_i)=
Q'_{ij}(\x_i)\,\dl_{ik}-r_{ik}\,q_{kj}\,.\label{plkq}\ee
Now multiplying by $d\x_k$ and summing over $k$ give (\ref{apde2}).

It follows from (\ref{plrho}) that $\pl_k\, \r_x=-R_x\,\dl_k\,\r$. Applying this to $A$  
gives $\pl_k\,Q'=-R_x\,\dl_k\,Q$, whose $i,j$ entry evaluated at $x=\x_i$
equals $-r_{xik}\,q_{kj}$. Hence
\be\pl_k\,q_{ij}'=\pl_k\,Q'_{ij}(\x_i)=-r_{xik}\,q_{kj}+\dl_{ik}\,Q''_{ij}(\x_i).
\label{plq'}\ee
Now we use Lemma~2 again. Applying both sides to $A$ and using the fact that 
$(D^2-M)A=0$ we obtain
\be Q''(x)-x\,Q(x)=R\,\dl\,Q'-R_y\,\dl\,Q.\label{Q''eq}\ee
Taking the $i,\,j$ entry and evaluating at $x=\x_i$ gives
\[Q_{ij}''(\x_i)-\x_i\,q_{ij}
=(rq'-r_{y}\,q)_{ij}.\]
Substituting this into (\ref{plq'}) we obtain
\[\pl_k\, q'_{ij}=-r_{xik}\,q_{kj}+\dl_{ik}\,[\x_i\,q_{ij}+
(rq'-r_{y}\,q)_{ij}].\]
Multipliying by $d\x_k$ and summing over $k$ give (\ref{apde4}).
 
To obtain the other equations, we point out that identities such as these occur in
dual pairs. Observe that the function $\ch_j(y)\,\r_{jk}(y,x)$ is equal
to $\ch_k(x)$ times $\tl\r_{kj}(x,y)$, where $\tl\r$ is the resolvent kernel
for the matrix kernel with $i,j$ entry $L_{ji}(x,y)\,\ch_j(y)$. Hence
$\tl Q_{jk}(x)$ is equal to $\ch_k(x)$ times the $Q_{kj}(x)$ associated with $L_{ji}$. The upshot 
is that for any formula involving $q$ or $\tl q$ there is another. 
We replace $q$ by $\tl q^t$ and  $\tl q$ with 
$q^t$. (If a formula involves $r$ we replace it by $r^t$ and subscripts $x$ and $y$ appearing in 
$r$ are interchanged.) In this way equations (\ref{apde3}) and (\ref{apde5})  are 
consequences of (\ref{apde2}) and (\ref{apde4}).
\sp

Let us derive the system of equations found in \cite{tw22}. We introduce the 
differential operator $\cD=\sum_k\pl_k$. The system of equations is
\begin{eqnarray}
{\cal {D}}^2\,q&=&\x\,q+2\,q\,\on\,\tl q\,q-2\,[\t,\,r]\,q,\label{eq1}\\
{\cal {D}}^2\,\tl q&=&\tl q\,\x+2\,\tl q\,q\,\on\,\tl q-2\,\tl q\,[\t,\,r],\label{eq2}\\
\cD r&=&-q\,\on\,\tl q+[\t,\,r].\label{eq3}
\end{eqnarray}

This can in fact be thought of as a system of ODEs since if we replace $\x_1,\cdots,\x_m$ by
$\x_1+\x,\cdots,\x_m+\x$
then $\cD=d/d\x$ and the $\x_j$ are parameters in the equations. 

Equation (\ref{eq3}) follows upon summing over $k$ the coefficients
of the $d\x_k$ in (\ref{apde1}) and using (\ref{plrsum}). Similarly
(\ref{apde2}) gives $\cD q=q'-r\,q,$ so
\be\cD^2\, q=\cD q'+(q\,\on\,\tl q-[\t,\,r])\,q-r\,(q'-r\,q).\label{D2q}\ee
Finally, (\ref{apde4}) gives
\[\cD q'=-(r_x+r_y)\,q+\x\,q+r\,q'.\]
Substituting this into (\ref{D2q}) and using (\ref{plrsum}) again give (\ref{eq1}). We derive
(\ref{eq2}) similarly.
\sp

When $m=1$ (\ref{eq1}) is the Painlev\'e~II equation $q''=\x q+2\,q^3$.
\sp

We now consider the more general case where each $X_k$ is a finite union of intervals,
\[X_k=(\x_{k1},\,\x_{k2})\cup(\x_{k3},\,\x_{k4})\cup\cd.\]
We write $\pl_{kw}$ for $\pl/\pl \x_{kw}$. We have
\be\pl_{kw}\,K=(-1)^w L\,\dl_{kw}(y),\label{plkwK}\ee
where $\dl_{kw}(y)$ is the $m\times m$ diagonal matrix all of whose entries are 0 except for the $k$th,
which equals $\dl(y-\x_{kw})$. It follows that 
\[\pl_{kw}\log\,\det(I-K)=-{\rm Tr}\;(I-K)\inv\,\pl_{kw}K=
(-1)^{w+1}R_{kk}(\x_{kw},\,\x_{kw}).\]

The various $\x_{kw}$ are the independent variables. (We shall systematically 
use $u,\ v$ and $w$ as indices to order the end-points of the intervals of $X_i,\ X_j$ and $X_k$,
respectively.) We now
define the matrix functions $r,\ q,\ \tl q,\ q'$ and $\tl q'$ of the $\x_{kw}$ by
\[r_{iu,\,jv}=R_{ij}(\x_{iu},\,\x_{jv}),\ \ \ q_{iu,\,j}=Q_{ij}(\x_{iu}),\ \ \ 
\tl q_{i,\,jv}=\tl Q_{ij}(\x_{jv}),\]
and
\[q'_{iu,\,j}=Q'_{ij}(\x_{iu}),\ \ \ \tl q'_{i,\,jv}=\tl Q'_{ij}(\x_{jv}).
\footnote{At points of discontinuity we always take limits from inside $X_k$.}
\]
These will be the unknown functions in our PDEs. We also define $r_x$ and $r_y$ by
\[r_{x,\,iu,\,jv}=R_{xij}(\x_{iu},\,\x_{jv}),\ \ \ r_{y,\,iu,\,jv}=R_{yij}(\x_{iu},\,\x_{jv}).\]

Observe that $r,\ r_x$ and $r_y$ are square 
matrices with rows and columns indexed by the end-points $kw$ of the $X_k$ while 
$q,\ q',\ \tl q$ and $\tl q'$ are rectangular matrices. Further notation is
\be\x=\dg(\x_{kw}),\ \ \ d\x={\rm diag}\,((-1)^{w+1}\,d\x_{kw}),\ \ \ 
\dh\x=\dg(d\x_{kw}),\ \ \ \dl=\sum_{k,w}(-1)^{w+1}\,\dl_{kw}.\label{notations}\ee
These are all square matrices but $\x,\ d\x$ and $\hat d\x$ are indexed by the end-points
of the $X_k$ while $\dl$ is $m\times m$.

With these notations our system of equations is
\begin{eqnarray}
dr&=&-r\,d\x\,r+\dh\x\,r_x+r_y\,\dh\x,\label{ade1}\\
dq&=&\dh\x\,q'-r\,d\x\,q,\label{ade2}\\
d\tl q&=&\tl q'\,\dh\x-\tl q\,d\x\,r,\label{ade3}\\
dq'&=&\dh\x\,\x\,q-(\,r_x\,d\x+d\x\,r_y)\,q+d\x\,r\,q',\label{ade4}\\
d\tl q'&=&\tl q\,\x\,\dh\x-\tl q\,(d\x\,r_y+r_x\,d\x)+\tl q'\,r\,d\x.\label{ade5}
\end{eqnarray}

As before the right sides involve the diagonal entries of $r_x+r_y$ and the off-diagonal entries of $r_x$ and
$r_y$, and we must show that these are known.

It is easy to see that Lemmas 1 and 2 still hold with the new definition of $\dl$.
Lemma 1 gives
\[ r_{x,\,iu,\,jv}+r_{y,\,iu,\,jv}=-\sum_{k,\ell}q_{iu,\,k}\,\tl q_{\ell,\,jv}+
\sum_{k,w}(-1)^w r_{iu,\,kw}\,r_{kw,\,jv}+(\t_i-\t_j)\,r_{iu,\,jv}.\]
In matrix terms,
\[r_x+r_y=-q\,\on\,\tl q-r\,s\,r+[\t,\,r],\]
where $s=\dg((-1)^{w+1})$. Thus $r_x+r_y$ is known. 

What remains is to show that $r_{x,iu,\,jv}$ and
$r_{y,iu,\,jv}$ are known when $iu\ne jv$. 

{}From (\ref{MRcom}) we have, using (\ref{DRcom}) again,
\[(x-y)\,R(x,y)=-Q'(x)\,\on\,\tl Q(y)+Q(x)\,\on\,\tl Q'(y)+(R_x+R_y)\dl\r-R\dl(\r_x+\r_y)+[\t,\,R_x-R_y],\]
so
\[[M,\,R]=-Q'(x)\,\on\,\tl Q(y)+Q(x)\,\on\,\tl Q'(y)-Q(x)\,\on\,\tl Q\dl\r(y)+R\dl Q(x)\,\on\,\tl Q(y)\]
\[+[\t,\,R]\,\dl \r-R\dl\,[\t,\,R]+[\t,\,R_x-R_y].\]
It follows, as before, that $r_{x,\,iu,\,jv}-r_{y,\,iu,\,jv}$ is known when $i\ne j$ and so 
also are $r_{x,\,iu,\,jv}$ and $r_{y,\,iu,\,jv}$ individually. It remains to determine these when $i=j$ but $u\ne v$.

To do this we use the identity $[DM,\,R]=D\,[M,\,R]+[D,\,R]\,M$ to compute
\pagebreak
\[[DM,\,R]=-Q''(x)\,\on\,\tl Q(y)+Q'(x)\,\on\,\tl Q'(y)-Q'(x)\,\on\,\tl Q\dl\r(y)+R_x\dl Q(x)\,\on\,\tl Q(y)\]
\[+[\t,\,R_x]\,\dl \r-R_x\dl\,[\t,\,R]+[\t,\,R_{xx}-R_{xy}]
-yQ(x)\,\on\,\tl Q(y)+yR\dl\r+y[\t,\,R].\]
Next we use (\ref{Q''eq}), which is the same here.
This gives an expression for $Q''(x)$ which we substitute into the first term above to obtain
\[[DM,\,R]=-(x+y)\,Q(x)\,\on\,\tl Q(y)+Q'(x)\,\on\,\tl Q'(y)-Q'(x)\,\on\,\tl Q\dl\r(y)
+(R_x+R_y)\dl Q(x)\,\on\,\tl Q(y)\]
\[+R\,\dl\,Q'(x)\,\on\,\tl Q(y)+[\t,\,R_x]\,\dl \r-R_x\dl\,[\t,\,R]+[\t,\,R_{xx}-R_{xy}]
+yR\dl\r+y[\t,\,R].\]

The left side equals $xR_x+yR_y+R$ and its $i,i$ entry evaluated at $(\x_{iu},\,\x_{iv})$ equals
$\x_{iu}\,r_{x,\,iu,\,iv}+\x_{iv}\,r_{y,\,iu,\,iv}+r_{iu,iv}$. If we can compute this sum 
then we know 
$r_{x,\,iu,\,iv}$
and $r_{y,\,iu,\,iv}$ individually since we know $r_{x,\,iu,\,iv}+r_{y,\,iu,\,iv}$ and 
$\x_{iu}\ne \x_{iv}$.
To see that the corresponding right side is computable observe that the term arising from 
$R_x+R_y$ is known because of 
Lemma~1, and the diagonal entries of $[\t,\,R_{xx}-R_{xy}]$ are zero. Everything else
is easily seen to be computable except possibly the terms arising from the sum
$[\t,\,R_x]\,\dl \r-R_x\dl\,[\t,\,R]$. Its $i,i$ entry equals two times
\[\t_i\sum_{k,w}(-1)^wR_{xik}(x,\,\x_{kw})\,R_{ki}(\x_{kw},\,y)-
\sum_{k,w}(-1)^wR_{xik}(x,\,\x_{kw})\,\t_k\,R_{ki}(\x_{kw},\,y).\]
The two summands corresponding to $k=i$ cancel. The remaining terms evaluated at 
$(\x_{iu},\,\x_{iv})$ 
involve $r_{kw,\,iv}$ and $r_{x,\,iu,\,kw}$ with $k\ne i$, all of which are known.

This completes the demonstration that all terms on the right sides of our equations are known.
This was the hard part. With 
(\ref{plrho}) replaced by 
\[\pl_{kw}\, \r=(-1)^w\,R\,\dl_{kw}\,\r,\]
the derivation of the equations proceeds exactly as before, and need not be repeated.
\sp

\noi{\bf Remark 1}. One might wonder whether the systems of equations (\ref{apde1})--(\ref{apde5}) and 
(\ref{ade1})--(\ref{ade5}) are integrable in the sense that one can 
derive from the
equations themselves that the differentials of the right sides are zero. Because of the 
complicated expressions for
$r_x$ and $r_y$ we have not attempted to show this in general. For equations
(\ref{apde1})--(\ref{apde5}), where we have relatively simple expressions for the right sides,
we verified that this is so when $m=2$ or 3.
\sp

\noi{\bf Remark 2}. We point out how little the equations 
depend on the operator $L$, as long as we still define $K=L\ch$ with
$\ch=\dg(\ch_{X_k})$. Equation (\ref{ade1}) holds for any integral operator $L$. So does 
(\ref{ade2}) if $q$ is defined as before in terms of $Q=(I-K)\inv\ph$, where $\ph$
can be any function whatsoever. Similarly for $\tl q$ and (\ref{ade3}). Similarly also
for the right hand sides of (\ref{ade4}) and (\ref{ade5}) except for the first terms
$\dh\x\,\x\,q$ and $\tl q\,\x\,\dh\x$. What does depend on the specifics of $L$ are the 
following:\sp\\ 
\noi(i) The expressions for $r_x$ and $r_y$ in terms of the unknowns. We do not see these
explicitly in
the equations. This is where the choice of $\ph$ arises.\sp\\
\noi(ii) The first terms on the right sides of (\ref{ade4}) and (\ref{ade5}), which 
arise from the computation of $Q''$. (See (\ref{Q''eq}).) \sp\\
All our systems will have the same 
form as these, most of the equations being universal, i.e., independent of the specific 
$L$ or $\ph$.\footnote
{This splitting into universal and nonuniversal equations was also a feature of \cite{tw5}.}
In most cases there will be two functions such as $\ph$. That will add to the number of
equations but not their complexity. The main difficulty in all cases will be (i).

\setcounter{equation}{0}\renewcommand{\theequation}{5.\arabic{equation}}
\begin{center}{\bf V. PDEs for the extended Hermite kernel}\end{center}

We modify (\ref{Lhermite}) by setting
\[L_{ij}(x,y)=\left\{\begin{array}{ll}\sum\limits_{k=0}^{n-1}e^{(k-n)\,(\t_i-\t_j)}\,
\ph_k(x)\,\ph_k(y)&{\rm if}\ i\ge j,
\\&\\-\sum\limits_{k=n}^\iy e^{(k-n)\,(\t_i-\t_j)}\,\ph_k(x)\,\ph_k(y)&{\rm if}\ i<j.
\end{array}\right.\]
The extra factors $e^{-n\,(\t_i-\t_j)}$ do not change the determinant.

Again we consider first the case where $X_k=(\x_k,\iy)$. We define $R$ and $\r$ as before, 
and again
\[\pl_k\log\,\det(I-K)=R_{kk}(\x_k,\,\x_k).\]
Now we shall have 
more unknown functions. We set
\[\ph=(2n)^{1/4}\,\ph_n,\ \ \ \ps=(2n)^{1/4}\,\ph_{n-1},\]
and define 
\[Q=\r\,\ph,\ \ \ P=\r\,\ps,\ \ \ \tl Q=\ph\ch\,\r,\ \ \ \tl P=\ps\ch\,\r.\]
Our unknowns will be, in addition to $r_{ij}=R_{ij}(\x_i,\,\x_j)$, the matrix functions 
$q,\ \tl q,\,p$ and $\tl p$ given by
\[q_{ij}=Q_{ij}(\x_i),\ \ \ \tl q_{ij}=\tl Q_{ij}(\x_j),\ \ \ p_{ij}=P_{ij}(\x_i),\ \ \ 
\tl p_{ij}=\tl P_{ij}(\x_j),\]
and
\[q'_{ij}=Q'_{ij}(\x_i),\ \ \ \tl q'_{ij}=\tl Q_{ij}(\x_j),\ \ \ p'_{ij}=P'_{ij}(\x_i),\ \ \ 
\tl p'_{ij}=\tl P'_{ij}(\x_j).\]
Again $\x$ denotes the matrix $\dg(\x_k)$ and $d\x$ denotes $\dg(d\x_k).$

With these notations our system of equations is
\begin{eqnarray}
dr&=&-r\,d\x\,r+d\x\,r_x+r_y\,d\x,\label{hpde1}\\
dq&=&d\x\,q'-r\,d\x\,q,\label{hpde2}\\
d\tl q&=&\tl q'\,d\x-\tl q\,d\x\,r,\label{hpde3}\\
dq'&=&d\x\,(\x^2-2n-1)\,q-(r_x\,d\x+d\x\,r_y)\,q+d\x\,r\,q',\label{hpde4}\\
d\tl q'&=&\tl q\,(\x^2-2n-1)\,d\x-\tl q\,(d\x\,r_y+r_x\,d\x)+\tl q'\,r\,d\x.\label{hpde5}\\
dp&=&d\x\,p'-r\,d\x\,p,\label{hpde6}\\
d\tl p&=&\tl p'\,d\x-\tl p\,d\x\,r,\label{hpde7}\\
dp'&=&d\x\,(\x^2-2n+1)\,p-(r_x\,d\x+d\x\,r_y)\,p+d\x\,r\,p',\label{hpde8}\\
d\tl p'&=&\tl p\,(\x^2-2n+1)\,d\x-\tl p\,(d\x\,r_y+r_x\,d\x)+\tl p'\,r\,d\x.\label{hpde9}
\end{eqnarray}

By Remark 2 and duality (each equation for $q$ or $p$ giving rise to one for $\tl q$
or $\tl p$) all we have to show is that the diagonal entries of 
$r_x+r_y$ and the off-diagonal entries of $r_x$ and $r_y$ are known (i.e., expressible
in terms of the unknowns) and to derive 
(\ref{hpde4}) and (\ref{hpde8}).

We begin by finding a substitute for Lemma 1. We write $D^{\pm}$ for $D\pm M$.
\sp

\noi{\bf Lemma 3}. We have
\be\Dp L_{ij}-e^{\t_i-\t_j}\,L_{ij}\Dp=-\ps(x)\,\ph(y),\ \ \ 
e^{\t_i-\t_j}\,\Dm L_{ij}-L_{ij}\Dm=-\ph(x)\,\ps(y).\label{DL}\ee

\noi{\bf Proof}. Let $J$ be the operator on $L^2({\bf R})$ with kernel 
\[J(x,y)=\sum_{k=0}^{n-1}\s^k\,\ph_k(x)\,\ph_k(y),\] 
and set $a_k=\sqrt{k/2}$. We have the formulas
\[x\ph_k=a_{k+1}\,\ph_{k+1}+a_k\,\ph_{k-1},\ \ \ \ph_k'=-a_{k+1}\,\ph_{k+1}+a_k\,\ph_{k-1}.\]
Therefore
\[(x+\pl_x)\,J(x,y)=2\,\sum_{k=0}^{n-1}\s^k\,a_k\,\ph_{k-1}(x)\,\ph_k(y),\]
\[(y-\pl_y)\,J(x,y)=2\,\sum_{k=0}^{n-1}\s^k\,a_{k+1}\,\ph_{k}(x)\,\ph_{k+1}(y).\]
This gives
\[[\,(x+\pl_x)-\s\,(y-\pl_y)\,]\,J(x,y)=-2\,\s^{n}\,a_{n}\,\ph_{n-1}(x)\,\ph_{n}(y).\]
If we take $\s=e^{\t_i-\t_j}$ and multiply by $e^{-n\,(\t_i-\t_j)}$ we obtain the first identity
of (\ref{DL})
when $i\ge j$. If $\s<1$ and one takes $n\to\iy$ in the last identity 
for $J$ one gets zero for the right 
sides. It follows that replacing $\sum_{k=0}^{n-1}$ by $-\sum_{k=n}^\iy$ in its definition  
does not change the right 
side. Thus we obtain the identity for $i<j$ as well. The second identity of (\ref{DL}) is 
obtained from the first by taking adjoints and using the fact that $L_{ij}$ is self-adjoint.
\sp

We can now find the analogue (actually, analogues) of Lemma~1. Observe that since
$\t=\dg(\t_i)$ we have $e^{\t}=\dg(e^{\t_i})$.

\noi{\bf Lemma 4}. We have
\be e^{-\t}\Dp\, R-R\,e^{-\t}\Dp=-P(x)\,e^{-\t}\,\on\,\tl Q(y)+R\dl e^{-\t}\r,\label{Rcom1}\ee
\be e^{\t}\Dm \,R-R\,e^{\t}\Dm=-Q(x)\,e^{\t}\,\on\,\tl P(y)+R\dl e^{\t}\r,\label{Rcom2}\ee
\sp 
\noi{\bf Proof}. If we multiply the relations (\ref{DL}) on the right by $\ch$ 
and use the fact $[D^{\pm},\,\ch]=\dl$ we obtain
\[e^{-\t}\Dp\, K-K\,e^{-\t}\Dp=-e^{-\t}\,\ps(x)\,\on\,\ch(y)\,\ph(y)+L\dl e^{-\t},\]
\[e^{\t}\Dm \,K-K\,e^{\t}\,\Dm=-e^{\t}\,\ph(x)\,\on\,\ch(y)\,\ps(y)+L\dl e^{\t}.\]
We replace $K$ on the left by $K-I$ and left- and right-multiply by $\r$, and the result follows.
(We used the fact that $e^{\pm\t}$ commutes with the matrix functions $\ph$ and $\ps$.)
\sp

If we take $i,\,j$ entries in (\ref{Rcom1}) and (\ref{Rcom2}) and set $x=\x_i,\ y=\x_j$ we obtain
\[e^{-\t} r_x+r_y\,e^{-\t}=-e^{-\t}\,\x\,r+r\,e^{-\t}\,\x-p e^{-\t}\,\on\, \tl q+r\,e^{-\t}\,r,\]
\be e^{\t}\,r_x+r_y\,e^{\t}=e^{\t}\,\x\,r-r\,e^{\t}\,\x-q e^{\t}\,\on\,\tl p+r\,e^{\t}\,r.\label{rxry}\ee
The right sides here are known. If we add and subtract these identities and take 
$i,\,j$ entries we obtain
\be 2\,(\cosh \t_i\,r_{xij}+\cosh\t_j\,r_{yij})=\cdots,\label{rxry1}\ee
\be 2\,(\sinh \t_i\,r_{xij}+\sinh\t_j\,r_{yij})=\cdots,\label{rxry2}\ee
where the dots on the right represent known quantities. The first relation with $j=i$ gives $r_{xii}+r_{yii}$.
If the two relations are thought of a system of equations for $r_{xij}$ and $r_{yij}$ the determinant
of the system is nonzero when $i\ne j$. Therefore we can solve for $r_{xij}$ and $r_{yij}$ individually then.
\sp

What remains is to derive (\ref{hpde4}) and
(\ref{hpde8}). For this we need the analogue of Lemma~2.
\sp

\noi{\bf Lemma 5}. We have 
\be [D^2-M^2,\,\r]=R\,\dl\,\r_x-R_y\,\dl\,\r,\label{D2M2com}\ee
where $R_y(x,y)$ is interpreted as not containing a delta-function summand.
\sp

The proof is analogous to that of Lemma~2. Here we use the fact that 
$D^2-M^2$ commutes with $L$, a consequence of the fact that each
$\ph_k$ is an eigenfunction of $D^2-M^2$. 

Since $\ph$ is an eigenfunction
of $D^2-M^2$ with eigenvalue $-2n-1$ and $\ps$ an eigenfunction with
eigenvalue $-2n+1$ applying both sides of (\ref{D2M2com}) to $\ph$ and to $\ps$ 
gives
\be Q''-x^2\,Q+(2n+1)\,Q=R\dl Q'-R_y\dl Q,\label{Qder}\ee
\be P''-x^2\,P+(2n-1)\,P=R\dl P'-R_y\dl P.\label{Pder}\ee

We have (\ref{plq'}) here just as before. Taking the $i,j$ entry in (\ref{Qder}) and
evaluating at $x=\x_i\,$ gives
\[Q_{ij}''(\x_i)-(\x_i^2-2n-1)\,q_{ij}
=(r\,q'-r_{y}\,q)_{ij}.\]
Substituting this into (\ref{plq'}) we obtain
\[\pl_k\, q'_{ij}=-r_{xik}\,q_{kj}+\dl_{ik}\,[(\x_i^2-2n-1)\,q_{ij}+
(r\,q'-r_{y}\,q)_{ij}],\]
which is (\ref{hpde4}). Equation (\ref{hpde8}) is established in exactly the same way using 
(\ref{Pder}).
\sp

We can also derive a system analogous to equations (\ref{eq1})--(\ref{eq3}):
\begin{eqnarray}
{\cal {D}}^2\,q&=&(\x^2-2n-1)\,q-2\,\cD r\cdot q,\label{heq1}\\
{\cal {D}}^2\,\tl q&=&\tl q\,(\x^2-2n-1)-2\,\,\tl q\cdot\cD r,\label{heq2}\\
{\cal {D}}^2\,p&=&(\x^2-2n+1)\,p-2\,\cD r\cdot p,\label{heq3}\\
{\cal {D}}^2\,\tl p&=&\tl p\,(\x^2-2n+1)-2\,\,\tl p\cdot\cD r,\label{heq4}\\
\cD r&=&-r^2+r_x+r_y.\label{heq5}
\end{eqnarray}

These equations are not as simple as (\ref{eq1})--(\ref{eq3}) since the expressions for the entries of
$r_x+r_y$ are messy. The last equation we already know. The other equations are
derived as for Airy: 
Summing the coefficients of $d\x_k$ in (\ref{hpde2}) gives $\cD q=q'-r\,q,$ so
\be\cD^2\, q=\cD q'-\cD\,r\cdot q-r\,(q'-r\,q).\label{hD2q}\ee
Similarly (\ref{hpde4}) gives
\[\cD q'=(\x^2-2n-1)\,q+r\,q'-(r_x+r_y)\,q.\]
Substituting this into (\ref{hD2q}) and using (\ref{heq5}) again give (\ref{heq1}). We derive
(\ref{heq3}) similarly, and (\ref{heq2}) and (\ref{heq4}) are obtained by duality.
\sp

In case $m=1$ (\ref{rxry}) gives $r_x+r_y=r^2-pq$, and our system of equations becomes
\[r'=-pq,\ \ \ q''=(\x^2-2n-1)\,q+2\,q^2\,p,\ \ \ p''=(\x^2-2n+1)\,p+2\,p^2\,q.\]
{}From the last two we find $(pq'-qp')'=pq''-qp''=-2pq$, and by the first equation this is $2r'$.
Thus $pq'-qp'=2r$. Using this, and successively computing $r'',\ r'''$ and $r''''$ using the
differentiation formulas, we arrive at
\[r''''=4(\x^2-2n)r''+4\x r'-12r'r''-4r=(4(\x^2-2n)r')'-4(\x r)'-6(r'\,^2)',\]
and so
\[r'''=4(\x^2-2n)r'-4\x r-6r'\,^2.\]
This  is the third-order equation found in \cite{tw5} which integrates to Painlev\'e~IV.
\sp

We turn to the more general case where each $X_k$ is a finite union of intervals, and will
again use the notations (\ref{notations}). The equations are
\begin{eqnarray}
dr&=&-r\,d\x\,r+\dh \x\,r_x+r_y\,\dh\x,\label{hde1}\\
dq&=&\dh\x\,q'-r\,d\x\,q,\label{hde2}\\
d\tl q&=&\tl q'\,\dh\x-\tl q\,d\x\,r,\label{hde3}\\
dq'&=&\dh\x\,(\x^2-2n-1)\,q-(r_x\,d\x+d\x\,r_y)\,q+d\x\,r\,q',\label{hde4}\\
d\tl q'&=&\tl q\,(\x^2-2n-1)\,\dh\x-\tl q\,(d\x\,r_y+r_x\,d\x)+\tl q'\,r\,d\x.\label{hde5}\\
dp&=&\dh\x\,p'-r\,d\x\,p,\label{hde6}\\
d\tl p&=&\tl p'\,\dh\x-\tl p\,d\x\,r,\label{hde7}\\
dp'&=&\dh\x\,(\x^2-2n+1)\,p-(r_x\,d\x+d\x\,r_y)\,p+d\x\,r\,p',\label{hde8}\\
d\tl p'&=&\tl p\,(\x^2-2n+1)\,\dh\x-\tl p\,(d\x\,r_y+r_x\,d\x)+\tl p'\,r\,d\x.\label{hde9}
\end{eqnarray}

Nothing is new here except to establish that the terms involving $r_x$ and $r_y$ on the 
right are known. As usual those that occur
are the diagonal entries of $r_x+r_y$ and the off-diagonal entries of $r_x$ and $r_y$.
In our case the terms $r_{xij}$ and $r_{yij}$ in (\ref{rxry1}) and (\ref{rxry2}) are replaced
by $r_{x,\,iu,\,jv}$ and $r_{y,\,iu,\,jv}$ and the relations show that these are known
when $i\ne j$ and that the $r_{x,\,iu,\,iv}+r_{y,\,iu,\,iv}$ are known. It remains
to show that $r_{x,\,iu,\,iv}$ and $r_{y,\,iu,\,iv}$ are known when $u\ne v$.

{}From (\ref{Rcom2}), which says
\[[e^{\t}\Dm, \,R]=-Q(x)\,e^{\t}\,\on\,\tl P(y)+R\dl e^{\t}\r,\]
we deduce 
\[[(e^{\t}\Dm)^2,\,R]=e^{\t}\Dm\,(-Q(x)\,e^{\t}\,\on\,\tl P(y)+R\dl e^{\t}\r)+
(-Q(x)\,e^{\t}\,\on\,\tl P(y)+R\dl e^{\t}\r)\,e^{\t}\Dm.\]
We use $S\equiv T$ for matrix functions $S$ and $T$ to denote that the differences
$S_{iu,iv}(\x_{iu},\,\x_{iv})-T_{iu,iv}(\x_{iu},\,\x_{iv})$ are known. 
If we keep in mind that $q,\ q',\ \tl p$ and $\tl p'$ are among our unknowns, we see
that it follows from the above, after multiplying by $e^{-2\t}$, that
\[[D^2-2MD,\,R]\equiv R_x\dl e^{\t}\r e^{-\t}-R\dl e^{\t} \r_y e^{-\t}.\]
If we subtract this from (\ref{D2M2com}) we obtain (since $[M^2,\,R]$ is known)
\be 2 [MD,\,R]\equiv (R\dl \r_x-R_x\dl e^{\t}\r e^{-\t})+(R\dl e^{\t} \r_y e^{-\t}-R_y\dl \r).
\label{MDhermitecom}\ee
Consider the first term on the right. Its $iu,\,iv$ entry evaluated at $(\x_{iu},\,\x_{iv})$
equals
\[\sum_{kw}r_{iu,\,kw}(-1)^{w+1}r_{x,\,kw,\,iv}-
\sum_{kw}r_{x,iu,\,kw}(-1)^{w+1}e^{\t_k}\,r_{\,kw,\,iv}\,e^{-\t_i}.\]
The terms of both sums corresponding to $k\ne i$ are known. So remaining as unknown is the
sum
\[\sum_{w}(-1)^{w+1}(r_{iu,\,iw}\,r_{x,\,iw,\,iv}-r_{x,iu,\,iw}\,r_{\,iw,\,iv}).\]
Analogously the second term on the right of (\ref{MDhermitecom}) is a known quantity plus
\[\sum_{w}(-1)^{w+1}(r_{iu,\,iw}\,r_{y,\,iw,\,iv}-r_{y,iu,\,iw}\,r_{\,iw,\,iv}).\]
Adding this to the last sum gives
\[\sum_{kw}(-1)^{w+1}[r_{iu,\,iw}\,(r_{x,\,iw,\,iv}+r_{y,\,iw,\,iv})-
(r_{x,\,iu,\,iw}+r_{y,iu,\,iw})\,r_{\,iw,\,iv})].\]
But this is known since, as we saw at the beginning, the $r_{x,\,iu,\,iv}+r_{y,iu,\,iv}$
are known.

We have shown that $[MD,\,R]\equiv$ a known matrix function. Its $iu,\,jv$ entry
evaluated at $(\x_{iu},\,\x_{iv})$ equals 
\[\x_{iu}\,r_{x,\,iu,\,jv}+\x_{iv}\,r_{y,\,iu,\,jv}+r_{y,\,iu,\,jv},\]
so $\x_{iu}\,r_{x,\,iu,\,jv}+\x_{iv}\,r_{y,\,iu,\,jv}$ is known. But so is
$r_{x,\,iu,\,jv}+r_{y,\,iu,\,jv}$. Therefore $r_{x,\,iu,\,jv}$ and $r_{y,\,iu,\,jv}$
are both known when $u\ne v$.
\pagebreak

\setcounter{equation}{0}\renewcommand{\theequation}{6.\arabic{equation}}
\begin{center}{\bf VI. PDEs for the extended sine kernel}\end{center}

If we make the substitutions $\t_i\to \t_i/2n,\ x\to x/\sqrt{2n},\ y\to y/\sqrt{2n}$ 
in the extended Hermite kernel (\ref{Lhermite}) and let $n\to\iy$  we obtain
the {\it extended sine kernel}
\[ L_{ij}(x,y)=\left\{\begin{array}{ll}{\displaystyle{\int_0^1}}e^{z^2\,(\t_i-\t_j)}\,
\cos z (x-y)\,dz&{\rm if}\ i\ge j,
\\&\\-{\displaystyle{\int_1^\iy}e^{z^2\,(\t_i-\t_j)}}\,\cos z (x-y)\,dz&{\rm if}\ i<j.
\end{array}\right.\]

Here we set 
\[\ph(x)=\sin x,\ \ \ps(x)=\cos x,\]
and then the other definitions are exactly as in Hermite with the above replacements. The unknowns now 
will be only $r,\ q,\ \tl q,\ p$ and $\tl p$ and the equations for general $X_k$ are
\begin{eqnarray}
dr&=&-r\,d\x\,r+\dh \x\,r_x+r_y\,\dh\x,\label{sde1}\\
dq&=&\dh\x\,(p+rsq)-r\,d\x\,q,\label{sde2}\\
d\tl q&=&(\tl p +\tl q s r)\,\dh\x-\tl q\,d\x\,r,\label{sde3}\\
dp&=&\dh\x\,(-q+rsp)-r\,d\x\,p,\label{sde4}\\
d\tl p&=&(-\tl q +\tl p s r)\,\dh\x-\tl p\,d\x\,r.\label{sde5}
\end{eqnarray}
(Recall that $s=\dg((-1)^{w+1})$\,.)

We know that equation (\ref{sde1}) is completely general, as are the equations
\begin{eqnarray*}
dq&=&\dh\x\,q'-r\,d\x\,q,\\
d\tl q&=&\tl q'\,\dh\x-\tl q\,d\x\,r,\\
dp&=&\dh\x\,p'-r\,d\x\,p,\\
d\tl p&=&\tl p'\,\dh\x-\tl p\,d\x\,r.
\end{eqnarray*}
To derive (\ref{sde2})--(\ref{sde5}) from these we establish
the formulas
\be q'=p+rsq,\ \ \ p'=-q+rsp,\ \ \ 
\tl q'=\tl p +\tl q s r,\ \ \ \tl p'=-\tl q +\tl p s r.\label{qp'}\ee

We have $[D,\,L]=0$, whence $[D,\,K]=L\dl$, whence
\be R_x+R_y=[D,\,\r]=R\,\dl\,\r.\label{DRcomsine}\ee
Applying (\ref{DRcomsine}) on the left to $\ph$ and $\ps$, using
$\ph'=\ps,\ \ps'=-\ph$, we obtain
\[Q'(x)-P(x)=(R\dl\,Q)(x),\ \ \ P'(x)+Q(x)=(R\dl\,P)(x).\]
Since  $q'_{iu,\,j}=Q'_{ij}(\x_{iu})$ and $p'_{iu,\,j}=P'_{ij}(\x_{iu})$, the first two
relations of
(\ref{qp'}) follow, and the others are analogous.

So equations (\ref{sde1})--(\ref{sde5}) hold, and it remains to deal with the entries of 
$r_x$ and $r_y$ appearing on the
right side of (\ref{sde1}). We have to show that the diagonal entries of
$r_x+r_y$ are known and that $r_{x,\,iu,\,jv}$ and $r_{y,\,iu,\,jv}$
are known when $iu\ne jv$.

It follows from (\ref{DRcomsine}) that $r_x+r_y=rsr$, and so all entries of the sum are known.

Next, for $i\ge j$ we have
\[2\,(\t_i-\t_j)\,L_{xij}=-2\,(\t_i-\t_j)\,\int_0^1
e^{z^2\,(\t_i-\t_j)}\,z\,\sin z (x-y)\,dz\]
\[=-e^{\t_i-\t_j}\,\sin (x-y)+(x-y)\,L_{ij}.\]
The same holds when $i<j$. Since $[\t D,\,L]_{ij}=\t_i\,L_{xij}+\t_j\,L_{yij}$ and
$L_y=-L_x$, this gives
\[2\,[\t D,\,L]=-e^{\t}\otimes e^{-\t}\,\sin(x-y)\,+[M,\,L],\]
where $e^{\t}\otimes e^{-\t}$ is the matrix with $i,j$ entry $e^{\t_i-\t_j}$. (This is not
a tensor product.) Hence
\[2\,[\t D,\,K]=-e^{\t}\otimes e^{-\t}\,\sin(x-y)\,\ch(y)+[M,\,K]+2\,K\dl\t.\]
Replacing $K$ by $K-I$ in the commutators and applying $\r$ left and right  
give
\[2\,[\t D,\,R]=P(x)\,e^{\t}\otimes e^{-\t}\,\tl Q(y)
-Q(x)\,e^{\t}\otimes e^{-\t}\,\tl P(y)+[M,\,R]+2\,R\dl\t\r.\]
The $i,j$ entry of the left side evaluated at $(\x_{iu},\,\x_{jv})$ equals twice
$\t_i\,r_{x,\,iu,\,jv}+\t_j\,r_{y,\,iu,\,jv}$\,, so these are known. We deduce, since
$r_{x,\,iu,\,jv}+r_{y,\,iu,\,jv}$ is known, that $r_{x,\,iu,\,jv}$ and $r_{y,\,iu,\,jv}$ are 
both known when $i\ne j$.
Just as before, the trickier part is to show that $r_{x,\,iu,\,iv}$ and $r_{y,\,iu,\,iv}$
are known when $u\ne v$.

We compute
\[[DM,\,R]=[D,\,R]\,M+D\,[M,\,R]\]
\be=R\dl\r\,y+2\,(\pl_x\,[\t D,\,R]-R_x\dl\t\r)-
P'(x)\,e^{\t}\otimes e^{-\t}\,\tl Q(y)
+Q'(x)\,e^{\t}\otimes e^{-\t}\,\tl P(y).\label{DMRcom}\ee
The $i,i$ entry of $[\t D,\,R]$ evaluated at $(\x_{iu},\,\x_{jv})$ equals
$\t_i\,(R_x+R_y)(\x_{iu},\,\x_{iv}).$ Hence, since $\pl_x\,(R_x+R_y)=R_x\dl\r$ by 
(\ref{DRcomsine}),
the $i,i$ entry of $\pl_x\,[\t D,\,R]$ evaluated at $(\x_{iu},\,\x_{jv})$ equals
$\t_i\,(R_x\dl\r)_{ii}(\x_{iu},\,\x_{iv})$. It follows that $i,i$ entry of 
$\pl_x\,[\t D,\,R]-R_x\dl\t\r$ evaluated at $(\x_{iu},\,\x_{jv})$ equals
\[\sum_{kw}(-1)^{w+1}(\t_i-\t_k)\,r_{x,\,iu,\,kw}\,r_{kw,\,iv}.\]
Since we need only sum over $k\ne i$ all these terms are known. So are the other terms of
(\ref{DMRcom}) evaluated at $(\x_{iu},\,\x_{jv})$.

We have shown that the $i,i$ entry of $[DM,\,R]$ evaluated at $(\x_{iu},\,\x_{jv})$ is known.
This equals $r_{iu,\,iv}+\x_{iu}\,r_{x,\,iu,\,jv}+\x_{iv}\,r_{y,\,iu,\,jv}$.
Thus $\x_{iu}\,r_{x,\,iu,\,jv}+\x_{iv}\,r_{y,\,iu,\,jv}$ is known. Since
$r_{x,\,iu,\,jv}+r_{y,\,iu,\,jv}$ is known and $\x_{iu}\ne\x_{iv}$ so also
are $r_{x,\,iu,\,jv}$ and $r_{y,\,iu,\,jv}$ known.
\sp

Let us see what these give in the case $m=1$ for a 
single interval $(-t,\,t)$. Here \linebreak $\x_1=-t,\ \x_2=t$. If we use the fact that $K(-x,\,-y)=
K(x,y)$ and the evenness of cosine and the oddness of sine we get $q_2=-q_1,\ p_2=p_1$
and if we use the fact that $R(x,\,y)=R(y,\,x)$ for $x,\,y\in(-t,\,t)$ we get 
$r_{ij}=r_{ji}$. 

We use
the notations $r=r_{11},\ \bar r=r_{12}$. If we observe that $d/dt=\pl_2
-\pl_1$ then (\ref{qp'}) gives 
\[{dq_1\ov dt}=-p_1-2\,\bar r\,q_1,\ \ \ {dp_1\ov dt}=q_1+2\,\bar r\,p_1,\]
and (\ref{sde1}) gives
\[{dr\ov dt}=r^2+{\bar r}^2-r_x-r_y\]
and the trivial relation $d{\bar r}/dt=-{\bar r}_x+{\bar r}_y$. 
The general relation $r_x+r_y=rsr$ gives in the present notation $r_x+r_y=r^2-{\bar r}^2$,
and so
\[{dr\ov dt}=2\,{\bar r}^2.\]
Finally (\ref{DMRcom}) gives
\[{\bar r}-t\,{\bar r}_x+t\,{\bar r}_y=-P_1'(-t)\,Q_2(t)+Q_i'(-t)\,P_2(t)
=-{d\ov dt}(Q_1(-t)\,P_1(-t)).\]
Thus 
\[{d\ov dt}(t\,{\bar r})=-{d\ov dt}(q_1\,p_1),\]
which gives
\[{\bar r}=-{q_1\,p_1\ov t}.\]
\sp

\setcounter{equation}{0}\renewcommand{\theequation}{7.\arabic{equation}}
\begin{center}{\bf VII. The Laguerre process}\end{center}

The Dyson process $\tau\rightarrow A(\tau)$ on the space of $p\times n$  complex matrices
(we assume $p\ge n$) is specified by its finite-dimensional distribution functions.
The probability measure on $A_k=A(\tau_k)\ (k=1,\ldots,m)$ is a normalization constant times
(\ref{Ldistr}), which may be written
\begin{eqnarray}
\prod_{j=1}^{m}\exp\left(-\left({1\over 1-q_j^2}+{q_{j+1}^2\over 1-q_{j+1}^2}\right)\,
\tr A_j^* A_j\right)\nonumber\\
\times \prod_{j=2}^m\exp\left({q_j\over 1-q_j^2}\,\tr(A_j^*A_{j-1}+{\rm hc})\right)
\, dA_1\cdots dA_{m}.\label{measure1} 
\end{eqnarray}
(Here ``hc'' is an abbreviation for ``Hermitian conjugate''.) We show how to derive
(\ref{Lprob}) from this.

Any complex matrix $p\times n$ complex 
matrix $A$ can, by the singular value
decomposition (SVD) theorem, be written as
\[ A= U D V^* \]
where $U$ is a $p\times p$ unitary matrix, $V$ is an $n\times n$ unitary matrix
and $D$ is a $p\times n$ matrix all of whose entries
are zero except for the diagonal consisting of the singular values
 of $A$. Thus we write each $A_j$ as
\[ A_j = U_j D_j V_j^* \]
with the goal of eventually integrating over the unitaries $U_j$ and $V_j$.
Of course,
\[ \tr(A_j^* A_j)=\tr(D_j^* D_j) = \sum_{j=1}^n \lambda_j \]
where $\lambda_j=d_{jj}^2$.

Let us examine one term
\[ \tr\left(A_j^* A_{j-1}+{\rm hc}\right) \]
appearing in the exponential of the second product in (\ref{measure1}).
Using the SVD representation we have terms
\[ \tr\left(V_j\,D_j^*\, U_j^* \,U_{j-1}\, D_{j-1}\, V_{j-1}^*+{\rm hc}\right) \]
The integrals over the unitary group (Haar measure) are both left- and right-invariant.
Thus in the $V_{j-1}$ integration we let
\[ V_{j-1}^*\rightarrow V_{j-1}^* V_j^* \]
so that the trace term becomes 
\[ \tr\left(D_j^*\, U_j^*\, U_{j-1}\, D_{j-1}\, V_{j-1}^*+{\rm hc}\right). \]
In the $U_{j-1}$ integration we  let $U_{j-1}\rightarrow U_{j} U_{j-1}$ and the trace becomes
\[ \tr\left(D_j^*\, U_{j-1}\, D_{j-1}\, V_{j-1}^*+{\rm hc}\right). \]
Thus, we have integrals of the form
\[ \int\int\exp\left({q_j\over 1-q_j^2}\,\tr(D_j^*\, U_{j-1}\, D_{j-1}
\, V_{j-1}^*+{\rm hc})\right)\,\,d\mu(V_{j-1})\, d\mu(U_{j-1})\]

Let $S$ denote an $n\times p$ complex matrix, $T$ a $p\times n$ complex matrix
and $U$ (resp.\ $V$) elements of the unitary group of $p\times p$  (resp.\
$n\times n$) matrices.  We assume $p\ge n$ and set $\alpha=p-n$.
We let $s_i$ resp.\ $t_i$ denote the eigenvalues
of $SS^*$ resp.\ $T^*T$. The
Harish-Chandra/Itzykson-Zuber integral 
for rectangular matrices (see, e.g., \cite{ZJZ}) is
\[
\int\int \exp\left(c\, \tr(S U T V^*+{\rm hc})\right)\,\, d\mu(U)\,d\mu(V)= 
{C_{p,n,c}\over \Delta(a) \Delta(b)}\, 
{\det\left(I_\alpha(2c\sqrt{a_i b_j}\,)
\right)\ov \prod_{i=1}^n (s_i \,t_i)^{\alpha/2}}.\]
Here $c$ can be any constant, $\al=p-n,\ I_\alpha$ is the modified
Bessel function and $C_{p.n,c}$ is a known constant.

When the $q_j=0$ the measure (\ref{measure1}) must reduce, after integration over
the unitary parts, to the well-known Laguerre measure.  It follows that (\ref{measure1})
becomes after integration over the unitary parts a normalization constant times
\pagebreak
\[\prod_{k=1}^{m}
e^{-\left({1\over 1-q_{k-1}^2}+{q_k^2\over 1-q_k^2}\right)\sum\limits_{i=1}^n\lambda_{ki}}\,
 \prod_{k=1}^{m-1}\det\left(I_\alpha\left({2 q_{k+1}\over 1-q_{k+1}^2}
\sqrt{\la_{i,\,k}\,\la_{j,k+1}}\right)\right)\]
\[\times\D(\la_{1})\, \D(\la_{m}) \,
\prod_{i=1}^n \la_{1i}^{\alpha/2}\,\prod_{i=1}^n \la_{mi}^{\alpha/2}\,\,
d\la_{11}\cdots d\la_{mn},\]
which is (\ref{Lprob}).
  
We shall now compute the extended kernel using the method of Section II. This density is not
quite of the form (\ref{density}) because of the last factors in the integrand here. Consequently in 
(\ref{ev}) there are extra factors $\la_i^{\al/2}$ and $\la_m^{\al/2}$, and so in the discussion that
follows $P_{1i}(\la)$ and $Q_{mi}(\la)$ are no longer polynomials of degree $i$ but $\la^{\al/2}$ times
polynomials of degree $i$.

We have now
\[ V_k(\lambda)=\left({1\over 1-q_{k-1}^2}+{q_k^2\over 1-q_k^2}\right) \lambda,
\>\>\> u_k(\lambda,\mu)=I_\alpha\left({2q_k\over 1-q_k^2} \sqrt{\lambda\mu}\right).\]
We introduce the Hille-Hardy kernel (the analogue of the Mehler kernel)
\[ K(q;\lambda,\mu)={q^{-\alpha}\over 1-q^2}\, e^{-{q^2\over 1-q^2}\,\lambda
-{1\over 1-q^2}\, \mu}\, \left({\mu\over \la}\right)^{\alpha/2}\, I_\alpha\left({2q\over 1-q^2}
\sqrt{\lambda\mu}\right)\]
which has the representation
\[ K(q;\lambda,\mu)=\sum_{i=0}^\iy q^{2i}\,p_i^\alpha(\lambda)\,p_i^\alpha(\mu)\, 
\mu^\alpha e^{-\mu}, \]
where $p_i^\alpha$ are the Laguerre polynomials $L_i^\al$, normalized.
It follows that
\be \int_0^\iy K(q;\lambda,\mu)\, p_i^\alpha(\mu)\,d\mu=q^{2i}\, p_i^\alpha(\la)\label{HH}\ee
and so again
$K(q)\ast K(q^\prime)=K(qq^\prime)$.

Now we may take in (\ref{factor}) 
\[E_{12}(\la_1,\la_2)=e^{-\la_1}\,(\la_1/\la_2)^{\al/2}\,K(q_1;\,\la_1,\la_2),\]
\[E_{k,k+1}(\la_k,\la_{k+1})=(\la_k/\la_{k+1})^{\al/2}\,K(q_k;\,\la_k,\la_{k+1}),\ \ \  (k>1),\]
and so
\[E_{1k}(\la,\mu)=e^{-\la}\,(\la/\mu)^{\al/2}\,K(q_1\cd q_{k-1};\,\la,\mu).\]
We deduce from (\ref{HH}) that
\[\int\int \,\la^{\al/2}\,p_i^\al(\la)\,E_{1m}(\la,\mu)\,\mu^{\al/2}\,p_j^\al(\mu)\,d\mu d\la\]
\[=(q_1\cd q_{m-1})^{2j}\int \la^\al\,p_i^\al(\la)\,e^{-\la}\,p_j^\al(\la)\,d\la=(q_1\cd q_{m-1})^{2j}\,\dl_{ij}.\]
Hence we may take
\[P_{1i}(\la_1)=\la_1^{\al/2}\,p_i^\al(\la_1),\ \ Q_{mj}(\la_m)=
(q_1\cd q_{m-1})^{-2j}\,\la_m^{\al/2}\,p_j^\al(\la_m).\]
We see that
\[P_{ki}(\mu)=\int \la^{\al/2}\,P_{1i}(\la)\,e^{-\la}\,(\la/\mu)^{\al/2}\,K(q_1\cd q_{k-1};\,\la,\mu)\,d\la\]
\[=(q_1\cd q_{k-1})^{2i}\,p_i^\al(\mu)\,\mu^{\al/2}\,e^{-\mu},\ \ \ (k>1),\]
\[Q_{kj}(\la)=\int (\la/\mu)^{\al/2}\,K(q_k\cd q_{m-1};\,\la,\mu)\,Q_{mj}(\mu)\,\mu^{\al/2}\,d\mu\]
\[=(q_1\cd q_{k-1})^{-2j}\,\la^{\al/2}\,p_j^\al(\la),\ \ \ (k>1),\]
\[Q_{1j}(\la)=\int e^{-\la}(\la/\mu)^{\al/2}\,K(q_1\cd q_{m-1};\,\la,\mu)\,Q_{mj}(\mu)d\mu
=e^{-\la}\,\la^{\al/2}\,p_j^\al(\la).\]

It follows that $H$ is the matrix with $k,\l$ entry
\[\sum_{j=0}^{n-1}\left({q_1\cd q_{\l-1}\ov q_1\cd q_{k-1}}\right)^{2j}\,p_j(\la)\,p_j(\mu)\]
left-multiplied by the matrix $\dg(\la^{\al/2}e^{-\la}\ \ \la^{\al/2}\ \cd\ \la^{\al/2})$ and right-multiplied
by the matrix \linebreak $\dg(\mu^{\al/2}\ \ \mu^{\al/2}e^{-\mu}\ \cd \ \mu^{\al/2}e^{-\mu})$. Similarly
$E$ is the strictly upper-triangular matrix with $k,\l$ entry 
$(\la/\mu)^{\al/2}\,K(q_k\cd q_{\l-1};\,\la,\mu)$
left-multiplied by the matrix $\dg(\la^{\al/2}e^{-\la}\ \ \la^{\al/2}\ \cd\ \la^{\al/2})$. 
Now we use the fact
that the determinant is unchanged if we multiply on the left by \linebreak
$\dg(e^{\la/2}\ \ e^{-\la/2}\ \cd\ e^{-\la/2})$ and on the right by 
$\dg(e^{-\mu/2}\ \ e^{\mu/2}\ \cd \ e^{\mu/2})$.

In this way we find the analogue of the kernel which was denoted by $\hat H-\hat E$ in Section III. 
It is now given by (\ref{Lhermite})
but with coefficients $e^{2k\,(\t_i-\t_j)}$ and with $\ph_k(x)$ equal to $x^{\al/2}\,e^{-x/2}\,p_k^\al(x)$.
This is the extended Laguerre kernel.\sp

\setcounter{equation}{0}\renewcommand{\theequation}{7.\arabic{equation}}
\begin{center}{\bf VIII. PDEs for the extended Bessel kernel}\end{center}

If we make the substitutions $\t_i\to\t_i/2n,\ x\to x^2/4n,\ y\to y^2/4n$ 
in the extended Laguerre kernel and then let $n\to\iy$ we obtain the {\it extended Bessel kernel}
\[ L_{ij}(x,y)=\left\{\begin{array}{ll} {\displaystyle{\int_0^1}} e^{z^2\,(\t_i-\t_j)/2}\,
\Ph_{\al}(xz)\,\Ph_{\al}(yz)\,dz&{\rm if}\ i\ge j,\\&\\ {\displaystyle -\int_1^\iy}
e^{z^2\,(\t_i-\t_j)/2}\,\Ph_{\al}(xz)\,\Ph_{\al}(yz)\,dz&{\rm if}\ i<j,
\end{array}\right.\]
where
\[\Ph_\al(z)=\sqrt z\,J_{\al}(z).\] 

Let us immediately explain the difficulty. In the previous cases we were able to find one
commutator for $L$ involving $D$ and another involving $D^2$, the latter arising from
the differential operator whose eigenfunctions appear in the integrand or summand of the expression for
the kernel. (For the extended Airy kernel these were given in Lemmas~1 and 2.) These enabled 
us to express $r_x$ and 
$r_y$ in terms of the unknown functions.

Here there does not seem to be
a commutator involving the first power of $D$. We are able to find two relations
involving the first power of $D$, but each involves both a commutator and an anticommutator.
Fortunately we are able to deduce from these relations three
commutator relations involving $D^2$, and these relations will enable us to show that
the derivatives of $r_x$ and $r_y$ are expressible in terms of $r_x$ and $r_y$ and
the other unknown functions. The upshot is that we are able to find a system of PDEs
in which $r_x$ and $r_y$ are now among the unknowns. Although the system of equations
seems no more complicated than those we have already derived (just larger) it is
actually much more so because of the expressions for
the derivatives of $r_x$ and $r_y$ in terms of the unknown functions.

To state the equations, we define $\ph$ and $\ps$ by
\[\ph=\Ph_{\al},\ \ \ \ \ps=\Ph_{\al+1}.\]
{}From these we define $q$ in the usual way. But now we set 
\[P=(I-K)\inv M\ps,\ \ \ \tl P=M\ps\ch\,(I-K)\inv,\]
and from these we define $p$ and $\tl p$ in the usual way. (The reason we do this is that 
eventually it is these $p$ and $\tl p$ which will arise in the expressions for the derivatives of $r_x$ and $r_y$.) 
With these notations our system of equations, in the general case where each $X_k$ is
a finite union of intervals, is\sp
\begin{eqnarray}
dr&=&-r\,d\x\,r+\dh \x\,r_{x}+r_{y}\,\dh\x,\label{bde1}\\
dr_x&=&-r_x\,d\x\,r+\dh \x\,r_{xx}+r_{xy}\,\dh\x,\label{bde2}\\
dr_y&=&-r\,d\x\,r_y+\dh \x\,r_{xy}+r_{yy}\,\dh\x,\label{bde3}\\ 
dq&=&\dh\x\,q'-r\,d\x\,q,\label{bde4}\\
d\tl q&=&\tl q'\,\dh\x-\tl q\,d\x\,r,\label{bde5}\\ 
dq'&=&\dh\x\,((\al^2-\textstyle{{1\ov4}})\,\x^{-2}\,q-q)-(r_x\,d\x+d\x\,r_y)\,q+d\x\,r\,q',\label{bde6}\\
d\tl q'&=&((\al^2-\textstyle{{1\ov4}})\,\tl q\,\x^{-2}-\tl q)\,\dh\x-\tl q\,(d\x\,r_y+r_x\,d\x)+\tl q'\,r\,d\x.\label{bde7}\\
dp&=&\dh\x\,p'-r\,d\x\,p,\label{bde8}\\
d\tl p&=&\tl p'\,\dh\x-\tl p\,d\x\,r,\label{bde9}\\
dp'&=&\dh\x\,((\al^2-\textstyle{{1\ov4}})\,\x^{-2}\,p+2q-p)-(r_x\,d\x+d\x\,r_y)\,p+d\x\,r\,p',\label{bde10}\\
d\tl p'&=&((\al^2-\textstyle{{1\ov4}})\,\tl p\,\x^{-2}+2\tl q-\tl p)\,\dh\x-\tl p\,(d\x\,r_y+r_x\,d\x)+\tl p'\,r\,d\x.\label{bde11}
\end{eqnarray}

Equations (\ref{bde2}) and (\ref{bde1}) are obtained in the same way as (\ref{ade1}). We have
\[\pl_{kw}r_{x,iu,jv}=\pl_{kw}R_{xij}(\x_{iu},\x_{jv})=
(-1)^w\,R_{xik}(\x_{iu},\x_{kw})\,R_{kj}(\x_{kw},\x_{jv})\]
\[+R_{xxij}(\x_{iu},\x_{jv})\dl_{iu,kw}+R_{xyij}(\x_{iu},\x_{jv})\dl_{jv,kw}.\]
This gives (\ref{bde2}) and (\ref{bde3}) is analogous.

So all the equations are universal except for (\ref{bde6}) and (\ref{bde10}) and their duals.
What we have to do is show that the diagonal entries of $r_{xx}+r_{xy}$ and $r_{xy}+r_{yy}$, 
and the off-diagonal entries of $r_{xx},\ r_{xy}$ and $r_{yy}$ are all known, and to
establish equations (\ref{bde6}) and (\ref{bde10}). 

To begin, we denote by $L^\pm$ the kernels where $\Ph_\al(xz)\,\Ph_\al(yz)$ in the integrand is 
replaced by
\[\Ph_\al(xz)\,\Ph_\al(yz)\pm \Ph_{\al+1}(xz)\,\Ph_{\al+1}(yz).\]

When $\al=-1/2,\ L^+$ is essentially the extended sine kernel and some of the formulas we derive here
will specialize to those obtained in Section~VI. We use the notations $\b={1\over2}+\alpha$
and
\[\z(x,y)=\ph(x)\,\ps(y)-\ps(x)\ph(y),\ \ \ \ \e(x,y)=\ph(x)\,\ps(y)+\ps(x)\,\ph(y),\]
\[\O=(e^{(\t_i-\t_j)/2}).\]

After integration by parts and some computation using the differentiation 
formulas
\[\Ph_\al'(z)=-\Ph_{\al+1}(z)+\b\, z\inv\Ph_\al(z),\ \ \ \
\Ph_{\al+1}'(z)=\Ph_\al(z)-\b\,z\inv\Ph_{\al+1}(z)\]
we find that
\[L^+_x={1\ov\t_i-\t_j}\O\,\z(x,y)+{\b\ov x}\,L^-+{1\ov\t_i-\t_j}(x-y)\,L^+,\]
\[L^+_y=-{1\ov\t_i-\t_j}\O\,\z(x,y)+{\b\ov y}\,L^-+{1\ov\t_i-\t_j}(y-x)\,L^+,\]
\[L^-_x=-{1\ov\t_i-\t_j}\O\,\e(x,y)+{\b\ov x}\,L^++{1\ov\t_i-\t_j}(x+y)\,L^-,\]
\[L^-_y=-{1\ov\t_i-\t_j}\O\,\e(x,y)+{\b\ov y}\,L^++{1\ov\t_i-\t_j}(x+y)\,L^-.\]
Here the $i,j$ entries of the matrices $L^\pm$ and $\O$ are to be understood. 

If we add the first two identities and subtract the last two we obtain we obtain the 
commutator-anticommutator pair 
\be[D,L^+]=\b\{M\inv,L^-\},\ \ \ \{D,L^-\}=\b[M\inv,L^+].\label{pair1}\ee

To obtain another pair, first multiply the first two identities by $\t_i-\t_j$
and subtract, getting
\[(\t_i-\t_j)\,(L_x^+-L_y^+)=2\,\O\,\z(x,y)+
\left({\b\ov x}-{\b\ov y}\right)(\t_i-\t_j)\,L^-
+2(x-y)\,L^+.\]
Using the first two identities again we can write the left side as
\pagebreak
\[\t_i\,L_x^++\t_j\,L_y^+-\t_i\left(-{1\ov\t_i-\t_j}\O\,\z(x,y)+{\b\ov y}\, L^-
+{1\over \t_i-\t_j}\,(y-x)\,L^+\right)\]
\[-\t_j\left({1\ov\t_i-\t_j}\O\,\z(x,y)+{\b\ov x}\,L^-+{1\over \t_i-\t_j}\,(x-y)\,L^+\right)\]
\[=\t_i\,L_x^++\t_j\,L_y^+-\left({\b\ov y}\t_i+{\b\ov x}\t_j\right)L^-+(x-y)\,L^+
+\O\,\z(x,y).\]
Thus
\[\t_i\,L_x^++\t_j\,L_y^+=\left({\b\ov x}\t_i+{\b\ov y}\t_j\right)L^-+(x-y)\,L^+
+\O\,\z(x,y).\]
In other words
\[[\t D-M,\,L^+]=\b\{M\inv\t,\,L^-\}+\O\,\z(x,y).\]

Next multiply the last two identities by $\t_i-\t_j$ and add, getting
\[(\t_i-\t_j)\,(L_x^-+L_y^-)=-2\,\O\,\e(x,y)
+\left({\b\ov x}+{\b\ov y}\right)(\t_i-\t_j)\,L^++2(x+y)\,L^-.\]
The left side may be rewritten
\[\t_i\,L_x^--\t_j\,L_y^-+\t_i\left(-{1\ov\t_i-\t_j}\O\,\e(x,y)+{\b\ov y}\, L^+
+{1\over \t_i-\t_j}\,(x+y)\,L^-\right)\]
\[-\t_j\left(-{1\ov\t_i-\t_j}\O\,\e(x,y)+{\b\ov x}\,L^++{1\over \t_i-\t_j}\,(x+y)\,L^-\right)\]
\[=\t_i\,L_x^--\t_j\,L_y^-+\left({\b\ov y}\t_i-{\b\ov x}\t_j\right)L^++(x+y)\,L^-
-\O\,\e(x,y).\]
Thus
\[\t_i\,L_x^--\t_j\,L_y^-=\left({\b\ov x}\t_i-{\b\ov y}\t_j\right)L^++(x+y)\,L^-
-\O\,\e(x,y).\]
In other words
\[\{\t D-M,\,L^-\}=\b[M\inv\t,\,L^+]-\O\,\e(x,y).\]

Thus we have our second commutator-anticommutator pair 
\be[\t D-M,\,L^+]=\b\{M\inv\t,\,L^-\}+\O\;\z(x,y),\label{pair2a}\ee
\be\{\t D-M,\,L^-\}=\b[M\inv\t,\,L^+]-\O\;\e(x,y).\label{pair2b}\ee

Now we have the following.
\sp

\noi{\bf Lemma}. Suppose $A$ and $B$ are such that
\[[A,\,L^+]=\{B,\,L^-\}+F,\ \ \ \ \{A,\,L^-\}=[B,\,L^+]+G.\]
Then 
\[[A^2-B^2,\,L^+]=[[A,\,B],\,L^-]+\{A,\,F\}+\{B,\,G\},\]
\[[A^2-B^2,\,L^-]=[[A,\,B],\,L^+]+[A,\,G]+[B,\,F].\]
\sp

\noi{\bf Proof}. We have 
$$[A^2,\,L^+]=\{A,\,[A,\,L^+]\}=\{A,\{B,\,L^-\}\}+\{A,\,F\}.$$
By the general identity
$$\{A,\{B,C\}\}=[[A,B],C]+\{B,\{A,C\}\}$$
the first term on the right side above may be written
$$[[A,B],L^-]+\{B,\{A,L^-\}\}=[[A,B],L^-]+\{B,[B,\,L^+]\}+\{B,\,G\}\]
\[=[[A,B],L^-]+[B^2,\,L^+]+\{B,\,G\}.$$
This establishes the first stated identity. 

For the second we write
\[[A^2,\,L^-]=[A,\,\{A,\,L^-\}]=[A,\,[B,\,L^+]]+[A,\,G].\]
By the general identity
\[[A,\,[B,\,C]]+[B,\,[C,\,A]]+[C,\,[A,\,B]]=0\]
the first term on the right side above may be written
\[-[B,\,[L^+,\,A]]-[L^+,\,[A,\,B]]=[B,\,\{B,\,L^-\}]+[B,\,F]-[L^+,\,[A,\,B]]\]
\[=[B^2,\,L^-]+[[A,\,B],\,L^+]+[B,\,F].\]
This gives the second identity.

We have obtained in (\ref{pair1}) and (\ref{pair2a})--(\ref{pair2b}) two
quadruples $(A_1,\,B_1,\,F_1,\,G_1)$ and \linebreak
$(A_2,\,B_2,\,F_2,\,G_2)$ satisfying the
hypothesis of the lemma. Each gives commutator relations involving $L^+$ and $L^-$.
However $(A_1+A_2,\,B_1+B_2,\,F_1+F_2,\,G_1+G_2)$ will also satisfy the 
hypothesis of the lemma and so gives commutator relations involving $L^+$ and $L^-$. If we
subtract from these the relations resulting from the other two we obtain 
\[[AA'+A'A-BB'-B'B,\,L^+]=[[A,\,B']+[A',\,B],\,L^-]+\{A,\,F'\}+\{A',\,F\}+\{B,\,G'\}+\{B',\,G\},\]
\[[AA'+A'A-BB'-B'B,\,L^-]=[[A,\,B']+[A',\,B],\,L^+]+[A,\,G']+[A',\,G]+[B,\,F']+[B',\,F].\]

So in the end we will obtain three pairs of commutator relations involving $L^+$ and $L^-$.
If we add the identites in each pair and divide by 2 we obtain three commutator identities
for $L$. For the explicit computations we have to keep in mind that all matrices and operators 
commute with $\ph$ and $\ps$,
and $D$ and $M$ commute with $\t$ and $\O$. We write down the results, sparing the reader 
the details:
\sp
\[[D^2+\b\,(1-\b)\,M^{-2},\,L]=0,\]
\[2\,[\t D^2-MD+\b(1-\b)\t\,M^{-2},\,L]\]
\[=\O\left(\ph\tn\ps\,(D-\b\,M\inv)-(D+\b\,M\inv)\,\ps\tn\ph\right),\]
\[[(\t\,D-M)^2+\b\,(1-\b)\,\t^2\,M^{-2},\,L]\]
\[=\O\,\ph\tn\ps\,(\t\,D-M-\b\,\t\,M\inv)-(\t\,D-M+\b\,\t\,M\inv)\,\O\,\ps\tn\ph.\]
The differentiation formula for $\Ph_{\al+1}$ is in our notation
$(D+\b\,M\inv)\,\ps=\ph$. Also, an operator acts on $\ph\tn\ps$ from the right
by applying its transpose to $\ps$. Using these facts we
see that the last two identities simplify to
\[[\t D^2-MD+\b(1-\b)\t\,M^{-2},\,L]=-2\,\O\,\ph\tn\ph,\]
\[[(\t\,D-M)^2+\b\,(1-\b)\,\t^2\,M^{-2},\,L]=-2\O\,\t\,\ph\tn\ph-\O(\ph\tn M\ps-M\ps\tn\ph).\]

The commutator identities for $L$ lead as before to commutator identities for
$K=L\ch$. They are
\[[D^2+\b\,(1-\b)\,M^{-2},\,K]=L\,(\dl D+D\dl),\]
\[[\t D^2-MD+\b(1-\b)\t\,M^{-2},\,K]=-2\,\O\,\ph\tn\ph\ch
+L\,\left(\t\,(\dl D+D\dl)-M\dl\right),\]
\[[(\t\,D-M)^2+\b\,(1-\b)\,\t^2\,M^{-2},\,K]=-2\O\,\t\,\ph\tn\ph\ch-
\O(\ph\tn M\ps\ch-M\ps\tn\ph\ch)\]
\[+L\,\left(\t^2\,(\dl D+D\dl)-2\t\,M\dl\right).\]

We are ready to apply $\r=(I-K)\inv$ to both sides. The only functions
that appear on the right sides are $\ph$ and $M\ps$, which is why we define
\[Q=(I-K)\inv \ph,\ \ \ \tl Q=\ph\,\ch\,(I-K)\inv,\ \ \ P=(I-K)\inv M\ps,\ \ \ 
\tl P=M\ps\,\ch\,(I-K)\inv.\]
Then we deduce
\be[D^2+\b\,(1-\b)\,M^{-2},\,R]=R\,\dl\,\r_x-R_y\,\dl\,\r,\label{bcom1}\ee
\[[\t D^2-MD+\b(1-\b)\t\,M^{-2},\,R]=-2\,Q(x)\,\O\,\tl Q(y)
+R\,\t\dl\,\r_x-R_y\,\t\dl\,\r-R\,\x\,\r,\]
\[[(\t\,D-M)^2+\b(1-\b)\t^2\,M^{-2},\,R]=-2Q(x)\,\O\t\,\tl Q(y)-
Q(x)\,\O\,\tl P(y)+P(x)\,\O\,\tl Q(y)\]
\[+R\,\t^2\dl\,\r_x-R_y\,\t^2\dl\,\r-R\,\t\x\,\r.\]

We now show that the diagonal entries of $r_{xx}+r_{xy}$ and $r_{xy}+r_{yy}$, 
and the off-diagonal entries of $r_{xx},\ r_{xy}$ and $r_{yy}$ can all be expressed
in terms of the unknowns.

We use the symbol $\equiv$ here to mean
that the difference of the quantities on its left and right is expressible in terms 
of $Q,\ P,\ \tl Q,\ \tl P$ and $R$, but no derivatives of these functions. The 
three commutator identities above yield in this notation the relations
\begin{eqnarray}R_{xx}-R_{yy}&\equiv& R\dl \r_x-R_y\dl\r,\label{1}\\
\t_i\,R_{xx}-\t_j\,R_{yy}&\equiv& x\,R_x+y\,R_y+R\,\t\dl \,\r_x-R_y\,\t\dl\,\r,\label{2}\\
\t_i^2\,R_{xx}-\t_j^2\,R_{yy}&\equiv &2\,\t\,x\,R_x+2\, y\,R_y\,\t+R\,\t^2\dl\, \r_x-
R_y\,\t^2\dl\,\r.
\label{3}\end{eqnarray}

Consider first the case $i\ne j$. It follows from any pair of the above equations 
(everything now is to be evaluated at $(\x_{iu},\x_{jv})$) that
both $R_{xx}$ and $R_{yy}$ are known. If we call the right sides above $A,\ B$ and $C$ then
\[\left|\begin{array}{ccc}1&1&A\\\t_i&\t_j&B\\\t_i^2&\t_j^2&C\end{array}\right|\equiv 0.\]
If we differentiate with respect to $x$ we deduce that the sum of all terms involving $R_{xy}$
is known. (Since our unknowns involved up to one derivative, this is why in our definition of $\equiv$ we required that no derivatives were
involved in the difference.) This sum is
\[-\t_i\t_j(\t_j-\t_i)\,R_{xy}\dl R-(\t_j^2-\t_i^2)\,(yR_{xy}-R_{xy}\t\dl R)
+(\t_j-\t_i)\,(2\t_j y\,R_{xy}-R_{xy}\t^2\dl R).\]
Dividing this by $\t_j-\t_i$, evaluating at $(\x_{iu},\x_{jv})$ and expanding we obtain
\[\sum_{k,w}(-1)^w(\t_i-\t_k)\,(\t_j-\t_k)\,r_{kw,jv}\,r_{xy,iu,kw}+
(\t_j-\t_i)\,\x_{jv}\,r_{xy,iu,jv}.\]
The terms involving $k=i$ vanish, so equating the above with the known quantity it is 
equal to gives a system of equations (with $iu$ fixed) for the $r_{xy,iu,kw}$ with $k\ne i$. The $jv,kw$
entry of the matrix for the system is
\[(-1)^w\,(\t_i-\t_k)\,(\t_j-\t_k)\,r_{kw,jv}+(\t_j-\t_i)\,\x_{jv}\,\dl_{jv,kw}.\]

The determinant of this matrix is a polynomial in the entries of $r$ and $\x$. (We think of the
$\t_j$ as fixed.) In the expansion of the determinant one summand is 
$\prod_{jv}(\t_j-\t_i)\,\x_{jv}$. Every other summand will contain at least one $r_{kw,jv}$ factor. If
we look at the series expansions for these other summands valid for small $\x_{jv}$ (coming from the series for 
the Bessel functions and
the Neumann series for the resolvent), every term will be a product of powers of the $\x_{jv}$ and 
have as coefficient a negative integral power of $\Gamma(\al)$
times a rational function of $\al$. It follows that in the series expansion of the determinant the coefficient
of $\prod_{jv}\x_{jv}$ is nonzero. Thus the determinant cannot be identically zero. 

We have shown that if $i\ne j$ then $r_{xy,iu,jv}$ is expressible in terms of the unknown functions. 
It remains to consider the cases where $i=j$, and we always evaluate at $(\x_{iu},\x_{iv})$.
In this case (\ref{1}) shows that $R_{xx}-R_{yy}$ is known. Subtracting $\t_i$ times
(\ref{1}) from (\ref{2}) gives
\[0\equiv x\,R_x+y\,R_y+R\t\dl R_x-R_y\t\dl R-\t_i(R\dl R_x-R_y\dl R).\]
All terms here involving $\dl$ are sums over $k$. The terms involving $k\ne i$, even after
taking $\pl_x$ or $\pl_y$, are known, as we have shown.  Those involving $k=i$ cancel,
just as before. 
Hence applying $\pl_x$ and $\pl_y$ to the above and evaluating at $(\x_{iu},\x_{iv})$ shows that
\[\x_{iu}\,r_{xx,iu,iv}+\x_{iv}\,r_{xy,iu,iv}\ \ {\rm and}\ \ 
\x_{iu}\,r_{xy,iu,iv}+\x_{iv}\,r_{yy,iu,iv}\]
are known. Taking $v=u$ shows that both $r_{xx,iu,iu}+r_{xy,iu,iu}$
and $r_{xy,iu,iu}+r_{yy,iu,iu}$ are known. If $u\ne v$, using the fact that 
$r_{yy,iu,iv}-r_{xx,iu,iv}$ is 
known we see also that $r_{xx,iu,iv}$ and $r_{xy,iu,iv}$
are individually known. 

All that we have left to show are (\ref{bde6}) and (\ref{bde10}). For these we use 
(\ref{bcom1}) (the analogue here of
Lemma~2) and the facts
\[(D^2+\b\,(1-\b)\,M^{-2})\,\ph=-\ph,\ \ \ (D^2+\b\,(1-\b)\,M^{-2})\,M\ps=2\ph-M\ps,\]
which follow from the differentiation formulas.
(The first is just the differential equation satisfied by $\Ph_\al$; the second is a miracle.)
We use these to compute $Q_{iu,j}''(\x_{iu})$ and $P_{iu,j}''(\x_{iu})$ as for previous equations.
Thus, for example, to obtain (\ref{bde6}) we replace the term $\dh\x\,(\x^2-2n-1)q$ in 
(\ref{hde4}) by $\dh\x\,(-\b\,(1-\b)\,\x^{-2}\,q-q)$ and to obtain (\ref{bde10}) we replace the term 
$\dh\x\,(\x^2-2n+1)p$ in (\ref{hde8})
by $\dh\x\,(-\b\,(1-\b)\,\x^{-2}\,p+2q-p).$ Any reader who has come this far
can easily supply the details.
\sp

\begin{center}{\bf Acknowledgments}\end{center}

We thank Kurt Johansson for sending us his unpublished notes on the extended
Hermite kernel.  
This work was supported by National Science Foundation under grants DMS-0304414 (first
author) and DMS-0243982 (second author).

\end{document}